\documentclass[12pt]{article}
\usepackage[dvips]{epsfig}
\usepackage{graphicx}
\usepackage{amsthm,amsmath,amssymb}
\usepackage{mathtools}
\usepackage{latexsym}
\usepackage{graphics}
\usepackage{url}
\usepackage{placeins}
\usepackage{color}
\usepackage{bm}
\usepackage{comment}
\usepackage{statex}
\usepackage{pgfplots}

\usepackage[
            bookmarksnumbered=true,
            bookmarksopen=true,
            colorlinks=true,
            citecolor=green,
            linkcolor=green,
            anchorcolor=red,
            urlcolor=blue
            ]{hyperref}

\definecolor{blue}{rgb}{0,0,0.9}
\definecolor{red}{rgb}{0.9,0,0}
\definecolor{green}{rgb}{0,0.50,0.10}
\definecolor{violet}{rgb}{0.5804,0.0000,0.8275}  
\newcommand{\blue}[1]{\begin{color}{blue}#1\end{color}}

\def\@themcountersep{}

\makeatletter
\newcommand{\labeltext}[2]{%
  \@bsphack
  \csname phantomsection\endcsname % in case hyperref is used
  \def\@currentlabel{#1}{\label{#2}}%
  \@esphack
}
\makeatother

\newtheorem{THEO}{Theorem}[section]
\newtheorem{ALGo}[THEO]{Algorithm}
\newtheorem{CONJ}[THEO]{Conjecture}
\newtheorem{COND}[THEO]{Condition}
\newtheorem{CORO}[THEO]{Corollary}
\newtheorem{DEFI}[THEO]{Definition}
\newtheorem{EXAMP}[THEO]{Example}
\newtheorem{FACT}[THEO]{Fact}
\newtheorem{HYPO}[THEO]{Hypothesis}
\newtheorem{LEMM}[THEO]{Lemma}
\newtheorem{PROB}[THEO]{Problem}
\newtheorem{PROP}[THEO]{Proposition}
\newtheorem{REMA}[THEO]{Remark}
\newcommand{\theo}{\begin{THEO}}
\newcommand{\algo}{\begin{ALGo} \rm}
\newcommand{\cond}{\begin{COND}}
\newcommand{\conj}{\begin{CONJ}}
\newcommand{\coro}{\begin{CORO}}
\newcommand{\defi}{\begin{DEFI} \rm}
\newcommand{\examp}{\begin{EXAMP} \rm}
\newcommand{\fact}{\begin{FACT}}
\newcommand{\hypo}{\begin{HYPO} \rm}
\newcommand{\lemm}{\begin{LEMM}}
\newcommand{\prob}{\begin{PROB} \rm}
\newcommand{\prop}{\begin{PROP}}
\newcommand{\rema}{\begin{REMA} \rm}
\newcommand{\etheo}{\end{THEO}}
\newcommand{\ealgo}{\end{ALGo}}
\newcommand{\econd}{\end{COND}}
\newcommand{\econj}{\end{CONJ}}
\newcommand{\ecoro}{\end{CORO}}
\newcommand{\edefi}{\end{DEFI}}
\newcommand{\eexamp}{\end{EXAMP}}
\newcommand{\efact}{\end{FACT}}
\newcommand{\ehypo}{\end{HYPO}}
\newcommand{\elemm}{\end{LEMM}}
\newcommand{\eprob}{\end{PROB}}
\newcommand{\eprop}{\end{PROP}}
\newcommand{\erema}{\end{REMA}}

\def\0{\mbox{\bf 0}}
\def\1{\mbox{\bf 1}}
\def\2{\mbox{\bf 2}}
\def\3{\mbox{\bf 3}}
\def\4{\mbox{\bf 4}}
\def\5{\mbox{\bf 5}}
\def\6{\mbox{\bf 6}}
\def\7{\mbox{\bf 7}}
\def\8{\mbox{\bf 8}}
\def\9{\mbox{\bf 9}}

\def\cc{\mbox{\boldmath $c$}}
\def\d{\mbox{\boldmath $d$}}

\def\t{\mbox{\boldmath $t$}}
\def\u{\mbox{\boldmath $u$}}
\def\v{\mbox{\boldmath $v$}}
\def\w{\mbox{\boldmath $w$}}
\def\x{\mbox{\boldmath $x$}}

\def\z{\mbox{\boldmath $z$}}
\def\A{\mbox{\boldmath $A$}}
\def\B{\mbox{\boldmath $B$}}
\def\C{\mbox{\boldmath $C$}}
\def\D{\mbox{\boldmath $D$}}

\def\F{\mbox{\boldmath $F$}}
\def\G{\mbox{\boldmath $G$}}
\def\H{\mbox{\boldmath $H$}}
\def\I{\mbox{\boldmath $I$}}

\def\L{\mbox{\boldmath $L$}}
\def\M{\mbox{\boldmath $M$}}

\def\O{\mbox{\boldmath $O$}}
\def\P{\mbox{\boldmath $P$}}
\def\Q{\mbox{\boldmath $Q$}}

\def\U{\mbox{\boldmath $U$}}
\def\V{\mbox{\boldmath $V$}}

\def\X{\mbox{\boldmath $X$}}
\def\Y{\mbox{\boldmath $Y$}}

\def\BC{\mbox{$\cal B$}}

\def\FC{\mbox{$\cal F$}}

\def\inprod#1#2{\langle#1, \, #2\rangle}
\def\Inprod#1#2{{\Large \langle}#1, \, #2{\Large \rangle}}

\def\Real{\mbox{$\mathbb{R}$}}
\def\Integer{\mbox{$\mathbb{Z}$}}

\def\s0{\mbox{\scriptsize \boldmath $0$}}

\def\wFC{\mbox{$\widehat{\FC}$}}

\def\Real{\mathbb{R}}

\def\coneC{\mathbb{C}}

\def\coneF{\mathbb{F}}
\def\spaceH{\mathbb{H}}

\def\coneK{\mathbb{K}}
\def\coneJ{\mathbb{J}}
\def\coneH{\mathbb{H}}

\def\spaceV{\mathbb{V}}
\def\SymMat{\mathbb{S}}

\def\Integer{\mathbb{Z}}
\def\coneG{\mbox{\boldmath $\Gamma$}}

\begin{document}

\title{
%\blue{
Exact SDP relaxations for a class of quadratic programs with finite and infinite quadratic
constraints 
%}
}

\author{
\normalsize
Naohiko Arima\thanks{Independent researcher %\red{Ph.D.  Policy and Planning Science, University of Tsukuba, March 2022 
	({\tt nao$\_$arima@me.com}).}, \and \normalsize
Sunyoung Kim\thanks{Department of Mathematics, Ewha W. University, 52 Ewhayeodae-gil, Sudaemoon-gu, Seoul 03760, Korea 
			({\tt skim@ewha.ac.kr}). 
			 The research was supported  by   NRF 2021-R1A2C1003810.}, \and \normalsize
Masakazu Kojima\thanks{Department of Industrial and Systems Engineering,
	Chuo University, Tokyo 192-0393, Japan ({\tt kojima@is.titech.ac.jp}).
 	}
}

% \date{\normalsize Revised, \today}
\date{\today}

\maketitle 

% 1. $f$         2. $\setf$   3.$\bvarphi$  4.

% \input abstract.tex
%!TEX root = ./main.tex

\begin{abstract}
\noindent
We investigate exact semidefinite programming (SDP) relaxations for
 the problem of minimizing a nonconvex quadratic objective function over a feasible region defined
by both finitely and infinitely many nonconvex quadratic inequality constraints (semi-infinite QCQPs).
Sufficient conditions for the exactness of SDP relaxations for QCQPs with finitely many constraints have been 
extensively studied, notably
by Argue et al. (MOR, 48:100-126, 2023), Arima et al. (SIOPT, 34:3194-3211, 2024), and
Joyce and Yang (MP, 205:539-558, 2024).
In this work, we present three new sufficient conditions that generalize the  existing conditions in these works 
for both finite and semi-infinite QCQPs.
Specifically, we establish  relationships among the proposed 
and existing  conditions,  and  prove that one of 
the proposed conditions is  the weakest among them, 
since it is implied by all the others.
Illustrative examples  are also provided to demonstrate the effectiveness of the proposed conditions 
in comparison to the existing ones.
\end{abstract}


\vspace{0.5cm}

%\red{

%

%\input sect1Introduction.tex
%!TEX root = ./main.tex

\section{Introduction}

\label{section:Introduction}

\label{section:generalCOP} 

We begin by considering a general conic optimization problem (COP). 
Let $\spaceV$ be a finite-dimensional vector space equipped 
with an inner product $\inprod{\A}{\B}$ for every $\A,\B\in \spaceV$. 
For every closed cone $\coneC \subseteq \spaceV$, $\Q \in \spaceV$ and 
$\H \in \spaceV$, COP$(\coneC,\Q,\H)$ denotes 
the problem of minimizing $\inprod{\Q}{\X}$ subject to 
$\X \in \coneC$ and $\inprod{\H}{\X} = 1$, {\it i.e.}, 
\begin{eqnarray*}
\eta(\coneC,\Q,\H) & = & \inf \left\{ \inprod{\Q}{\X} : \X \in \coneC, \ \inprod{\H}{\X} = 1 \right\},
\end{eqnarray*}
where $\coneC \subseteq \spaceV$ is a cone if $\lambda \X \in \coneC$ holds for every 
$\X \in \coneC$ and $\lambda \geq 0$.
We note that a cone $\coneC$ is not 
necessarily convex. When $\Q, \H \in \spaceV$ are unspecified and arbitrary, we 
often denote COP$(\coneC,\Q,\H)$ and $\eta(\coneC,\Q,\H)$ as COP$(\coneC)$ and 
$\eta(\coneC)$, respectively. 
If COP$(\coneC)$ is infeasible, we assume that $\eta(\coneC) = +\infty$. 

\medskip

%%%%%
Let $\coneK$ be a closed nonconvex cone in $\spaceV$. Let  $\mbox{co}\coneK$ and 
$\overline{\mbox{co}}\coneK$ denote the convex hull of $\coneK$ and its closure, respectively.
For every closed convex cone $\coneJ \subseteq \mbox{co}\coneK$, $\Q \in \spaceV$ and 
$\H \in \spaceV$, we consider a nonconvex conic optimization problem 
COP$(\coneK \cap\coneJ)$:  
\begin{eqnarray*}
\eta(\coneK\cap\coneJ) 
&= &\inf\left\{ \inprod{\Q}{\X} : \X \in \coneK\cap\coneJ, \ \inprod{\H}{\X} = 1 \right\} 
\end{eqnarray*}
and its convex  relaxation  COP$(\coneJ)$:
 \begin{eqnarray*}
\eta(\coneJ) &= &\inf\left\{ \inprod{\Q}{\X} : \X \in \coneJ, \ \inprod{\H}{\X} = 1 \right\}.
\end{eqnarray*}
Obviously, $\eta(\coneJ) \leq \eta(\coneK\cap\coneJ)$ holds.
If $\eta(\coneJ) = \eta(\coneK\cap\coneJ)$, 
we say that 
COP($\coneK\cap\coneJ$) and its convex relaxation COP($\coneJ$)
are  {\em equivalent}, or that 
COP($\coneJ$)  is {\em an exact convex relaxation of} 
COP$(\coneK\cap\coneJ)$. 

\medskip

The above framework of a nonconvex conic optimization problem 
COP$(\coneK \cap\coneJ)$ and its convex relaxation COP($\coneJ$) 
was originally introduced in \cite{KIM2020} for a unified geometrical analysis on 
the completely positive programming (CPP) reformulation of quadratically constrained 
quadratic problems (QCQPs) and their extension to polynomial optimization 
problems (POPs). They introduced 
\begin{eqnarray*}
\wFC(\coneK) & = & \mbox{the family of closed convex cones $\coneJ \subseteq
\mbox{co}\coneK$
such that co$(\coneK\cap\coneJ) = \coneJ$},  
\end{eqnarray*}
and established the following results.
%%%%%%%%%%
\theo \label{theorem:propFC}
Let $\Q, \ \H \in \spaceV$, and $\coneK \subseteq \spaceV$ be a closed cone.
\vspace{-2mm}
\begin{description}
\item{(i) } Assume that $\coneJ \in \wFC(\coneK)$. Then 
\begin{eqnarray}
& & -\infty < \eta(\coneJ,\Q,\H) \ 
 \mbox{if and only if} \  -\infty < \eta(\coneJ,\Q,\H) = \eta(\coneK\cap\coneJ,\Q,\H).
 \label{eq:wFC}
 \end{eqnarray}
\vspace{-8mm}
\item{(ii) } If $\coneJ$ is a face of co$\coneK$,  then $\coneJ \in \wFC(\coneK)$. 
\vspace{-2mm}
\item{(iii) } $\coneJ' \in \wFC(\coneK)$ for every face $\coneJ'$ of 
$\coneJ \in \wFC(\coneK)$.
\vspace{-2mm}
\item{(iv) } Assume that $\H \in \mbox{int co}\coneK$ (the interior of co$\coneK$). 
Then $\coneJ \in \widehat{\FC}(\coneK)$ if and only if 
$\eta(\coneJ \cap \coneK,\Q,\H) = \eta(\coneJ,\Q,\H)$ for every $\Q \in \SymMat^{n}$.
\vspace{-1mm}
\end{description} 
\etheo
%%%%%%%%%%
\proof{
For assertion (i), we refer to \cite[Theorem 3.1]{KIM2020} (see also 
\cite[Corollary 2.2]{ARIMA2023});
 for (ii) and (iii), to \cite[Lemma 2.1]{KIM2020};  
 for (iv) 
 to \cite[Theorem 1.2]{ARIMA2023}. 
} 

\medskip

\noindent
In particular, Kim et al. \cite{KIM2020} applied Theorem~\ref{theorem:propFC} (i) and (ii) to 
equivalent CPP relaxation of a class of QCQPs 
in binary and nonnegative 
variables, which includes Burer's class of QCQPs \cite{BURER2009},  
and its extension to POPs.

\subsection{A quadratically constrained quadratic program (QCQP) and its 
equivalent semidefinite programming (SDP) relaxation}

\label{section:QCQPSDP}

In this paper, we focus on the case $\spaceV = \SymMat^n$, 
the linear space of $n \times n$ symmetric matrices, 
%This paper deals with a special case: $\spaceV = \SymMat^n$ 
%(the linear space of $n \times n$ symmetric matrices).
 and $\coneK = 
\coneG^n \equiv \{ \x\x^T : \x \in \Real^n\}$,  where $\Real^n  $ is the 
$n$-dimensional Euclidean space of column vectors $\x = (x_1,\ldots,x_n)$ 
and $\x^T$ denotes the row vector obtained by transposing $\x \in \Real^n $. 
In this case, \vspace{-2mm}
\begin{itemize}
\item co$\coneK = \mbox{co}\coneG^n = \SymMat^n_+$ (the convex cone 
of $n \times n$ positive semidefinite matrices),\vspace{-2mm}
\item  COP($\coneG^n\cap\coneJ$) and its convex relaxation COP($\coneJ$) 
correspond to a (geometric form of) QCQP and its semidefinite 
programming (SDP) relaxation,
\vspace{-2mm}
\item each closed convex cone 
$\coneJ \in \wFC(\coneG^n)$ is characterized as 
$\coneJ = \mbox{co}\left\{\x\x^T : \x \in \Real^n, \ \x\x^T \in \coneJ\right\}$, 
and is called {\em rank-one-generated (ROG)} in the literature 
\cite{ARGUE2023,BLEKHERMAN2017,HILDEBRAND2016}.\vspace{-1mm}
\end{itemize} 

\medskip

Argue et al. \cite{ARGUE2023}  demonstrated independently 
from \cite{KIM2020}, that the ROG property is crucial for
ensuring the equivalence relation~\eqref{eq:wFC} in case $\coneK = \coneG^n$  
\cite[Lemma 19]{ARGUE2023}. 
They provided a thorough study of characterizations of ROG cones and 
established several necessary and sufficient conditions. In particular, condition~(II) below 
characterizes when a closed convex cone $\coneJ \subseteq \SymMat^n_+$,
defined by linear matrix inequalities, has the ROG property. 
In \cite{ARIMA2023},  
Arima et al. also presented a sufficient condition, 
condition (I) below for $\coneJ \in \wFC(\coneG^n)$. 

\medskip

A closed convex cone $\coneJ \subseteq \SymMat^n_+$ is often represented 
using linear matrix inequalities.  For every $\B\in \SymMat^n$,  let 
\begin{eqnarray*}
& & \coneJ_+(\B), \ \coneJ_0(\B) \ \mbox{or  } \coneJ_-(\B) = \left\{ \X \in \SymMat^n_+ : \inprod{\B}{\X} \geq, \ = \mbox{or } \leq 0, \ 
\mbox{respectively}
\right\}, 
\end{eqnarray*}
and  $\coneJ_+(\BC) = \left\{ \X \in \SymMat^n_+ : \inprod{\B}{\X} \geq 0 \ (\B \in \BC) \right\}$  for every $\BC \subseteq \SymMat^n$. 
Since $\coneJ \subseteq \SymMat^{n}_+$ is a closed convex cone, 
$\coneJ$ can be represented as 
the intersection of (possibly infinitely many) half spaces and $\SymMat^{n}_+$ 
such that $\coneJ =\coneJ_+(\BC)$ for some $\BC \subseteq \SymMat^{n}$.
We should mention that there are many choices for such a $\BC$.  
For example, we can take $\BC = \coneJ^*$, where 
$\coneJ^* = \{ \Y \in \SymMat^n :  \inprod{\X}{\Y} \geq 0 \ \mbox{for every } \X \in \coneJ\}$ 
(the dual of $\coneJ$). 
This trivial choice of $\BC$, however, involves many redundant matrices to represent~$\coneJ$. 

\medskip

Now we represent a QCQP and its SDP relaxation as 
COP($\coneG^n\cap\coneJ_+(\BC)$) 
\begin{eqnarray}
\eta(\coneG^n\cap\coneJ_+(\BC)) & = & 
 \inf\left\{ \inprod{\Q}{\x\x^T} : 
\begin{array}{l}
\x \in \Real^n, \ \inprod{\B}{\x\x^T} \geq 0 \ (\B \in \BC),\\
\inprod{\H}{\x\x^T} = 1
\end{array}
\right\} \label{eq:QCQPBC0} \\[3pt]
& = & 
\inf\left\{ \inprod{\Q}{\X} : 
\begin{array}{l}
\X \in \SymMat^n, \X \in \coneG^n, \\
\X \in\coneJ_+(\B) \ (\B \in \BC), \ 
\inprod{\H}{\X} = 1
\end{array}
\right\}, \label{eq:QCQPBC} 
\end{eqnarray}
and COP($\coneJ_+(\BC)$) 
\begin{eqnarray}
\eta(\coneJ_+(\BC)) & = & \inf\left\{ \inprod{\Q}{\X} : 
\begin{array}{l}
\X \in \SymMat^n_+, \\
\X\in\coneJ_+(\B) \ (\B \in \BC), \
 \inprod{\H}{\X} = 1 
\end{array}
 \right\}, \label{eq:SDPBC} 
\end{eqnarray}
respectively. 
We note that $\x\in\Real^n$ is the variable in~\eqref{eq:QCQPBC0} while 
$\X\in\SymMat^n$ is the variable in~\eqref{eq:QCQPBC}. 
The process of embedding QCQP~\eqref{eq:QCQPBC0}, defined in $\x\in\Real^n$, 
into QCQP ~\eqref{eq:QCQPBC} 
is often called {\em a lifting into the matrix space $\SymMat^n$}. 
The following results are known:
%%%%%%%%%%
\theo \label{theorem:existing0}
Let $\BC \subseteq \SymMat^n$. If one of the following conditions (I) and (II)  holds, then 
$\coneJ_+(\BC) \in \wFC(\coneG^n)$. 
In particular, $\coneJ_+(\B) \in \wFC(\coneG^n)$, $\coneJ_-(\B) = \coneJ_+(-\B) \in 
\wFC(\coneG^n)$ and $\coneJ_0(\B) = \coneJ_+(\{\B,-\B\})\in\wFC(\coneG^n)$ 
for every $\B \in \SymMat^n$. 
\vspace{-1mm}
\etheo
%%%%%%%%%%
\begin{description}
\item{(I) } $\BC$ is finite. 
$\coneJ_0(\B) \subseteq \coneJ_+(\A)$ for every $\A, \ \B\in \BC$. 
\cite[Theorem 4.1]{ARIMA2023}. 
\vspace{-2mm}
\item{(II)} $\BC$ is finite. 
For every distinct $\A, \ \B\in \BC$, there exists a nonzero 
$(\alpha,\beta)  \in \Real^2$ such that $\alpha\A + \beta\B \in \SymMat^n_+$ 
\cite[Proposition 1]{ARGUE2023}.
\vspace{-2mm}
\end{description}

\medskip

For QCQP examples that satisfy condition~(I) in \cite{JEYAKUMAR2014,POLYAK1998,STERN1995,YE2003}, 
we refer the reader to \cite{ARIMA2023}. For fundamental properties of ROG cones and their applications to equivalent SDP relaxations of QCQPs,
see \cite{ARGUE2023,KARZAN2021}.

\medskip

\subsection{Non-intersecting quadratic constraint conditions}

\label{section:non-intersection}

%\indent
Conditions (I) and (II) are defined in $\SymMat^n$, the space of the variable matrix $\X$ of 
COP($\coneJ_+(\BC)$). 
For practical applications, however,
it is more convenient to provide a direct characterization of QCQP~\eqref{eq:QCQPBC0}, 
which involves $\BC$ satisfying condition (B),  in the space $\Real^n$ of 
its variable vector $\x$. 
We examine the special case where 
$\H = \H^1 \equiv \mbox{diag}(0,\ldots,0,1) \in \SymMat^n_+$ (the $n \times n$ 
diagonal matrix with 
diagonal elements $0,\ldots,0,1$). 
Letting  
\begin{eqnarray*}
\spaceH^1 = \{ \X \in \SymMat^n : \inprod{\H^1}{\X} = 1 \} =  \{ \X \in \SymMat^n : X_{nn} = 1 \}, 
\end{eqnarray*}
we can rewrite COP$(\coneG^n\cap\coneJ_+(\BC),\Q,\H^1)$ and 
COP$(\coneJ_+(\BC),\Q,\H^1)$ as 
\begin{eqnarray*}
\eta(\coneG^n\cap\coneJ_+(\BC),\Q,\H^1) & = & \inf 
\left\{ \inprod{\Q}{\X} : \X \in \coneG^n\cap\coneJ_+(\BC)\cap \spaceH^1\right\} \\[3pt]
& = & \inf 
\left\{ \inprod{\Q}{\X} : \X\in \overline{\mbox{co}}(\coneG^n\cap\coneJ_+(\BC)\cap \spaceH^1)\right\}
\\[3pt]
& & \mbox{(since the objective function $\inprod{\Q}{\X}$ is linear in $\X \in \SymMat^n$)}
\end{eqnarray*}
and 
\begin{eqnarray*}
\eta(\coneJ_+(\BC),\Q,\H^1) & = & \inf 
\left\{ \inprod{\Q}{\X} : \X \in \coneJ_+(\BC)\cap \spaceH^1\right\}, 
\end{eqnarray*}
respectively. 
Therefore, $\overline{\mbox{co}}(\coneG^n\cap\coneJ_+(\BC)\cap \spaceH^1) = 
\coneJ_+(\BC)\cap \spaceH^1$ serves as a sufficient condition for the equivalence of 
COP$(\coneG^n\cap\coneJ_+(\BC),\Q,\H^1)$ and COP$(\coneJ_+(\BC),\Q,\H^1)$, 
which will be utilized in Theorem~\ref{theorem:JOYCE2024} below.

\medskip

For the subsequent discussion, 
COP$(\coneG^n\cap\coneJ_+(\BC),\Q,\H^1)$ will be expressed in an alternative form.  
Define  
\begin{eqnarray*}
q(\u,\B) & = & 
\inprod{\B}{{\scriptsize \begin{pmatrix} \u \\ 1 \end{pmatrix}}
{\scriptsize \begin{pmatrix} \u \\ 1 \end{pmatrix}}^T}
\mbox{ for every } (\u,\B) \in \Real^{n-1}\times\SymMat^n.  
\end{eqnarray*}
If  $\B = {\scriptsize \begin{pmatrix} \C & \cc^T \\ \cc & \gamma \end{pmatrix}}$,
where  $\C \in \SymMat^{n-1}$, $\cc \in \Real^{n-1}$ and $\gamma \in \Real$, then $q(\u,\Q)$ is 
a quadratic function of  the form $\u^T\C\u + 2\cc^T\u + \gamma$ in $\u \in \Real^{n-1}$. 
For every $\B\in \SymMat^n$ and $\BC \subseteq \SymMat^n$, define 
\begin{eqnarray*}
& & 
\B_\geq, \ \B_= \ \mbox{or } \B_\leq  =  
\left\{ \u \in \Real^{n-1} :  q(\u,\B) \geq 0,  \ = \ \mbox{or } \leq 0, \ \mbox{respectively}
\right\}, \\
& & 
\BC_\geq = \bigcap_{\B \in \BC} \B_\geq = 
\left\{ \u \in \Real^{n-1} :  q(\u,\B) \geq 0 \ (\B \in \BC)\right\}.
\end{eqnarray*}
 Now, we rewrite QCQP \eqref{eq:QCQPBC0}
with $\H = \H^1 = \mbox{diag}(0,\ldots,0,1) \in \SymMat^n_+$ as 
\begin{eqnarray}
\eta(\coneG^n\cap\coneJ_+(\BC),\Q,\H^1) 
& = & \inf\left\{q(\u,\Q) : \u \in \Real^{n-1}, \ q(\u,\B) \geq 0 \ (\B \in \BC)\right\} \nonumber \\
& = &  \inf\{q(\u,\Q) : \u \in \BC_\geq \}. \label{eq:QCQPBCFH}
\end{eqnarray}

\medskip

The following result follows from \cite[Corollary 3]{JOYCE2024}, 
specialized to COP$(\coneG^n\cap\coneJ_+(\BC),\Q,\H^1)$ and 
COP$(\coneJ_+(\BC),\Q,\H^1)$.  
%%%%%%%%%%
\theo \label{theorem:JOYCE2024} 
Let $\BC \subseteq \SymMat^n$.  Assume that condition (III) holds. 
Then 
\begin{eqnarray}
\coneJ_+(\BC)\cap\spaceH^1=\overline{\mbox{co}}(\coneG^n\cap\coneJ_+(\BC)\cap\spaceH^1),
\label{eq:JOYCE2024a}
\\ [3pt]
\eta(\coneJ_+(\BC),\Q,\H^1) =  \eta(\coneG^n\cap\coneJ_+(\BC),\Q,\H^1).
\label{eq:JOYCE2024b}
\end{eqnarray}
%hold. 
(Moreover, \eqref{eq:JOYCE2024a} immediately implies \eqref{eq:JOYCE2024b}.)
%(Recall that~\eqref{eq:JOYCE2024a} implies~\eqref{eq:JOYCE2024b}).  
\etheo  
%%%%%%%%%%
\begin{description}
\item{(III) } $\BC$ is finite.  
$q(\cdot,\B) : \Real^{n-1} \rightarrow \Real$ % $(\B \in \BC)$ 
is not affine 
({\it i.e.}, if $q(\cdot,\B)$ is described as $q(\u,\B) = \u^T\D\u + 2\d^T\u + \delta$ 
for every $\u \in \Real^{n-1}$, then $\D \not= \O$) $(\B \in \BC)$ 
and 
\begin{eqnarray}
\B_=  \subseteq \BC_\geq \ \mbox{ for every } \B \in \BC.
\label{eq:condOYCE2024}
\end{eqnarray}
\end{description}

The condition~\eqref{eq:condOYCE2024} is 
commonly referred to as 
%frequently called 
{\em the non-intersecting quadratic constraint condition} (NIQCC) 
\cite{BURER2015,JOYCE2024,YANG2018}. 
Early works examined special instances of NIQCC for relatively simple QCQPs, 
most notably those arising from the generalized trust region subproblem (TRS)
%Early studies considered special cases of NIQCC for 
%simple QCQPs, particularly those arising from the generalized trust region subproblem (TRS) 
\cite{BECK2006,BURER2015,CONN2000,JEYAKUMAR2014,PONG2014}. 
Extending the generalized TRS framework,
%As an extension of the generalized TRS, 
%the works in
the authors in \cite{POLYAK1998,YE2003} studied % investigated  
 QCQPs 
 of the form $\inf\{q_0(\u) : -1 \leq q_1(\u) \leq 1\}$, 
where $q_0, \ q_1$ are  quadratic functions in $\u \in \Real^{n-1}$. 
This problem can be reformulated as 
a QCQP satisfying NIQCC. 
Similarly, 
a quadratic program with non-intersecting ellipsoidal hollows \cite{YANG2018} provides another extension of the generalized TRS. Thus, 
condition \eqref{eq:condOYCE2024} may be 
viewed as a unified formulation of NIQCC-type assumptions that arise in these classes of QCQPs. 
%regarded as a unified representation or generalization of NIQCC-type assumptions that appear in these QCQPs.
 %When compared with condition (I),  
 In comparison with condition (I), 
 \eqref{eq:condOYCE2024} 
 represents  %can be viewed as 
a {\em non-homogenized form of NIQCC, whereas 
condition (I) corresponds to } a %epresents a
{\em homogenized form of NIQCC}.

\subsection{Summary of main results}

The main purpose of the paper is to propose new conditions that generalize previously mentioned  conditions (I), (II) and (III)
and to investigate their relationships in detail. 
To describe the new conditions,  
we introduce a face $\coneF$ of $\SymMat^n_+$ that contains 
$\coneJ_+(\BC)$. 
Clearly, $\coneJ_+(\BC) = \coneJ_+(\BC)\cap\coneF$. 
To ensure that the conditions are nontrivial and meaningful, 
we implicitly assume that $\coneF$ is a proper face  of $\SymMat^n_+$.   
If a linear equality $\G{\scriptsize \begin{pmatrix}\u\\1\end{pmatrix}} = 0$ is given 
for some $r \times n$ full row-rank matrix $\G$, 
it can be incorporated into QCQP~\eqref{eq:QCQPBCFH} 
 as a quadratic constraint $q(\u,\B) \geq 0$, where $\B = -\G^T\G$. In this case, $\coneJ_+(\BC)$ 
 is contained in the face 
 $\coneF = \coneJ_+(-\G^T\G) = \coneJ_0(\G^T\G)$. (Note that $\G^T\G\in\SymMat^n_+$). 
Moreover, even if none of $\coneJ_+(\B)$ 
 $(\B \in \BC)$ is a face of $\SymMat^n_+$, 
  $\coneJ_+(\BC)$
 may still be contained in a proper face of
  $\SymMat^n_+$, as demonstrated in Section~\ref{section:exampleFBmain1}. 
Therefore, studying such cases is both theoretically and practically important.

\medskip

As a weaker variant of condition (I), we propose the following condition:\vspace{-2mm}
\begin{description}
\item{(B) } 
$\coneJ_0(\B) \cap \coneF \subseteq \coneJ_+(\A) \cap \coneF$ 
or $\coneJ_+(\A) \cap \coneF \subseteq \coneJ_+(\B) \cap \coneF$ holds for every 
$\A, \B \in \BC$. \vspace{-2mm}
\end{description}
%%%%%%%%%%
\theo \label{theorem:main1} 
Let $\BC \subseteq \SymMat^n$ and $\coneF$ be a face of  $\SymMat^n_+$ that includes 
$\coneJ_+(\BC)$. Assume that condition (B) holds.
Then $\coneJ_+(\BC) \in \wFC(\coneG^n)$. 
\etheo
%%%%%%%%%%
We emphasize that an  infinite set $\BC$ is allowed in condition (B); hence  QCQP~\eqref{eq:QCQPBC} 
may be  a semi-infinite program \cite{Hettich1993,REEMTSEN2023}. 
The inclusion 
$\coneJ_+(\A) \cap \coneF \subseteq \coneJ_+(\B) \cap \coneF$ 
in condition (B) for distinct $\A, \B\in \BC$ 
implies that $\B$ is redundant to describe $\coneJ_+(\BC)\cap\coneF$, {\it i.e.}, 
$\coneJ_+(\BC)\cap\coneF = \coneJ_+(\BC\backslash\{\B\})\cap\coneF$, 
but  $\coneJ_+(\BC)$ can be a proper subset of $\coneJ_+(\BC\backslash\B))$, 
as shown  in Section~\ref{section:exampleFBmain1}. 
Theorem~\ref{theorem:main1} strengthens \cite[Theorem 3.1]{ARIMA2023}  
where the authors assume that $\BC$ is finite and fix 
$\coneF = \SymMat^n_+$ as in condition (I). 

\medskip

It is well-known that every face of $\SymMat^n_+$ is isomorphic to $\SymMat^r_+$ for some 
$r \in \{0,\ldots,n\}$ \cite{PATAKI2000,PATAKI2013}. 
For proofs of Theorems~\ref{theorem:main1}, \ref{theorem:main2} 
and~\ref{theorem:JOYCE2024ext},  we apply facial reduction from $\coneF$ onto 
$\SymMat^r_+$. (See \cite{BORWEIN1981,BORWEIN1981a,WAKI2013} 
for numerical methods for facial  reduction).
Let $\Phi$ be a linear isomorphism from $\coneF \supseteq \coneJ_+(\BC)$ 
onto $\SymMat^r_+$, and 
$\Phi^* : \SymMat^n \rightarrow \SymMat^r$ the adjoint map with respect to 
$\Phi$.  
Then we can prove that 
$\Phi(\coneJ_0(\B)\cap\coneF) = \coneJ_0(\Phi^*(\B))$ and  
$\Phi(\coneJ_+(\B)\cap\coneF) = \coneJ_+(\Phi^*(\B))$ for every $\B\in \SymMat^n$,  and that 
$\Phi(\coneJ_+(\BC))  = \Phi(\coneJ_+(\BC)\cap\coneF) = \coneJ_+(\Phi^*(\BC))$. 
Thus condition (B) is equivalent to\vspace{-2mm}
\begin{description}
\item{($\widetilde{B}$) } 
$\coneJ_0(\widetilde{\B}) \subseteq \coneJ_+(\widetilde{\A})$ 
or $\coneJ_+(\widetilde{\A})  \subseteq \coneJ_+(\widetilde{\B})$ holds for every 
$\widetilde{\A}, \widetilde{\B} \in \Phi^*(\BC)$.\vspace{-2mm}
\end{description}
Also, $\coneJ_+(\BC) \in \wFC(\coneG^n)$ holds if and only if  
$\coneJ_+(\Phi^*(\BC)) = \Phi(\coneJ_+(\BC)) \in \Phi(\wFC(\coneG^n)) = \wFC(\coneG^r)$. 
As a result, when $\BC$ is finite, we can reduce Theorem~\ref{theorem:main1} to 
Theorem~\ref{theorem:existing0} as will be shown in Section~\ref{section:ProofMain1Finite}. 
We note that  the isomorphism $\Phi$ and the \blue{adjoint} $\Phi^*$ do not appear explicitly  
in the descriptions of the proposed conditions 
(B), (C) and~(D), although they play an essential role in the proofs of our main results, 
including Theorems~\ref{theorem:main1}, \ref{theorem:main2} 
and~\ref{theorem:JOYCE2024ext}. It should be emphasized that 
conditions (B), (C) and~(D) are invariant under the isomorphism $\Phi$, while 
conditions (I), (II) and (III) are not. 
Further details will be provided in Section~\ref{section:Main1}. 
The equivalence of $\coneJ_0(\BC) \in \wFC(\coneG^n)$ and $\coneJ_0(\Phi^*(\BC)) \in \wFC(\coneG^r)$, which was obtained similarly 
in \cite[Observation 1]{ARGUE2023}, 
was utilized for the proof of [1, Proposition 1] in a different context, where 
$\coneJ_0(\BC) = \{ \X \in \SymMat^n_+ : \inprod{\B}{\X} = 0 \ (\B \in \BC)\}$. But 
the condition (II) derived there is not invariant under $\Phi$.

\medskip

We also propose non-homogenized NIQCCs
(C) and (D) below. Let 
\begin{eqnarray*}
L(\coneF)= \{ \x \in \Real^n : \x\x^T \in \coneF \},\ 
L_1(\coneF)  = \{ \u \in \Real^{n-1} : 
{\scriptsize \begin{pmatrix} \u \\ 1\end{pmatrix}} \in L(\coneF)\}. 
\end{eqnarray*}
Since $\coneF$ is a face of $\SymMat^n_+$, there exists an $\F \in \SymMat^n_+$ 
such that $\coneF = \{ \X \in \SymMat^n_+: \inprod{\F}{\X} = 0\}$. 
Hence $L(\coneF)= \{ \x \in \Real^n : \inprod{\F}{\x\x^T} = 0 \} 
= \{ \x \in \Real^n : \F\x = \0 \}$, which implies that $L(\coneF)$ 
is a linear subspace of $\Real^{n}$, and $L_1(\coneF) = \{ \u \in \Real^{n-1} : \F {\scriptsize \begin{pmatrix} \u \\ 1\end{pmatrix}} = \0 \}$ an affine subspace of $\Real^{n-1}$. 
It follows from $\coneJ_+(\BC) \subseteq \coneF$ 
 that $\BC_\geq \subseteq L_1(\coneF)$; hence $\BC_\geq = \BC_\geq \cap L_1(\coneF)$.
Specifically, if $\coneF = \SymMat^n_+$ then 
$L(\coneF)= \Real^n$ and $L_1(\coneF) = \Real^{n-1}$.\vspace{-2mm}

\medskip

%%%%%%%%%%
\theo \label{theorem:main2} Let $\BC \subseteq \SymMat^n$ and $\coneF$ be a face of 
$\SymMat^n_+$ that includes $\coneJ_+(\BC)$. 
Assume that the following condition (C) is satisfied. Then condition (B) holds; hence $\coneJ_+(\BC) \in \wFC(\coneG^n)$.
\etheo
%%%%%%%%%%
\begin{description}
\item{(C) }  
$\BC_\geq \cap L_1(\coneF) \not= \emptyset$, and 
\begin{eqnarray*}
\emptyset \not= \B_\leq \cap L_1(\coneF) \subseteq \A_\geq\cap L_1(\coneF) \ \mbox{or }
\A_\geq \cap L_1(\coneF) \subseteq \B_\geq\cap L_1(\coneF) \ 
\mbox{for every } \A, \ \B  \in \BC. 
\end{eqnarray*}
\vspace{-10mm}
\end{description}

%%%%%%%%%%
\theo \label{theorem:JOYCE2024ext} 
Let $\BC \subseteq \SymMat^n$ and $\coneF$ be a face of 
$\SymMat^n_+$ that includes $\coneJ_+(\BC)$. 
Assume that the following condition (D) is satisfied. 
Then~\eqref{eq:JOYCE2024a} and~\eqref{eq:JOYCE2024b} hold.  
\etheo
%%%%%%%%%%
\begin{description}
\item{(D) }  
$\BC$ is finite,
$q(\cdot,\B) : \Real^{n-1} \rightarrow \Real$ $(\B \in \BC)$ is not affine 
on $L_1(\coneF)$, and 
 \begin{eqnarray}
\B_= \cap L_1(\coneF) \subseteq \A_\geq \cap L_1(\coneF) \ 
\mbox{ for every } \A, \ \B \in \BC. \label{eq:condJOYCE2024ext} 
 \end{eqnarray}
\end{description}
If we take $\coneF = \SymMat^n_+$ then condition (D) coincides with condition (III). 
Thus Theorem \ref{theorem:JOYCE2024ext}   is a generalization of Theorem~\ref{theorem:JOYCE2024}. 
Neither of conditions (C) and (D) implies the other;  
the inclusion relation $\B_\leq\cap L_1(\coneF) \subseteq \A_\geq\cap L_1(\coneF)$ in
(C) implies the inclusion relation $\B_=\cap L_1(\coneF) \subseteq \A_\geq\cap L_1(\coneF)$  
in (D), 
while (C) allows $q(\cdot,\B)$  for some $\B \in \BC$ to be affine 
on $L_1(\coneF)$. Figure 1 shows three cases (a), (b) and (c) for conditions~(C) and~(D), 
where condition (C) is satisfied in all cases but condition (D) is satisfied only 
in cases (a).  
Case (b) involves an affine function $q(\cdot,\B^5) = u_2 + 1$. 
In case (c), $\B^1_= \not\subseteq \B^4_\geq$, $\B^2_= \not\subseteq \B^4_\geq$ and 
$\{\B^1,\B^2,\B^4\}_\geq = \{\B^4\}_\geq$ hold. This case illustrates that 
any redundant $\B$ should be removed from $\BC$ for condition  
(D) to be more effective.  In conditions (C), 
`or $\A_\geq \cap L_1(\coneF) \subseteq \B_\geq\cap L_1(\coneF)$' is added to 
adapt such cases.

\medskip

\begin{figure}[t!]  \vspace{-20mm}  
\begin{center}
\includegraphics[height=85mm]{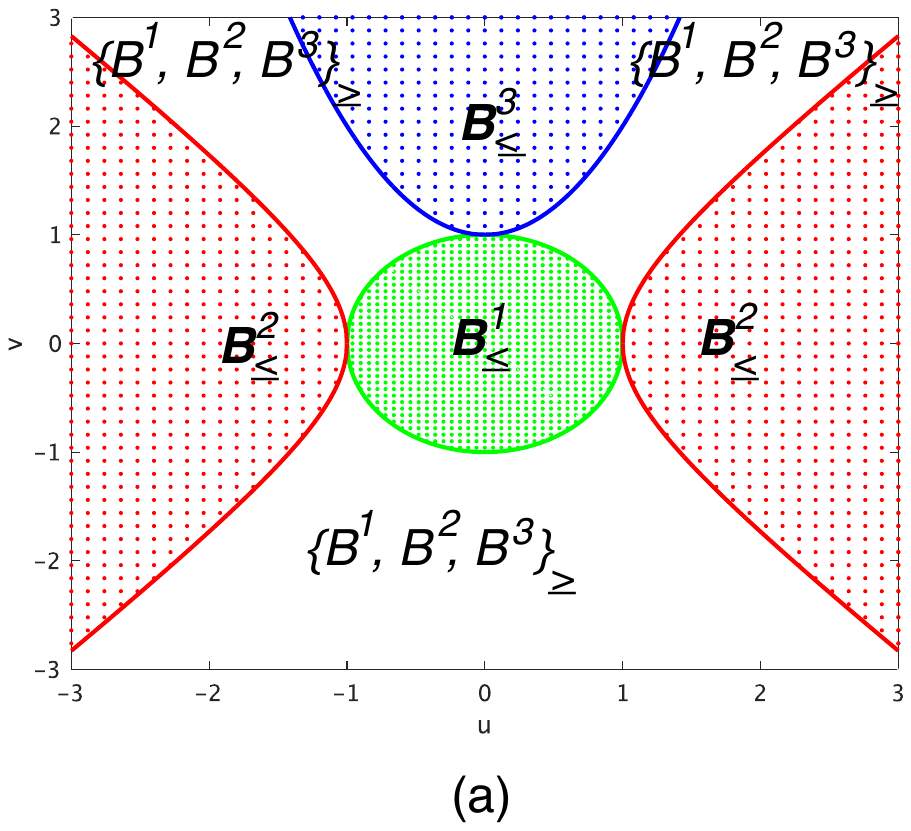} \hspace{-14mm}
\includegraphics[height=85mm]{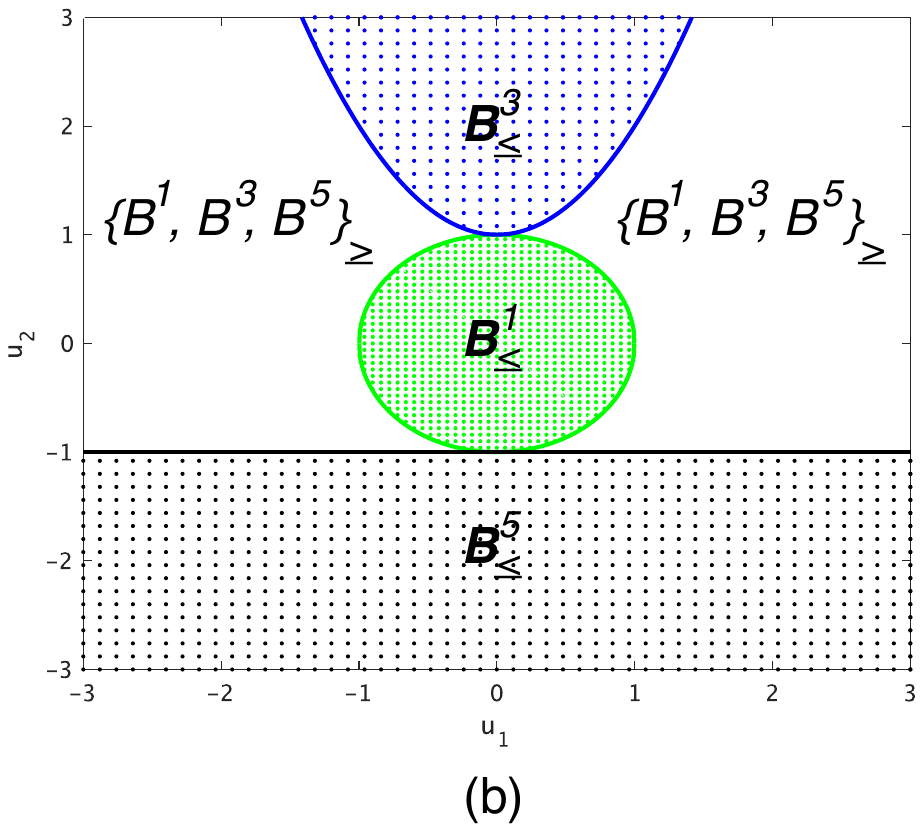}\hspace{-14mm}
\includegraphics[height=85mm]{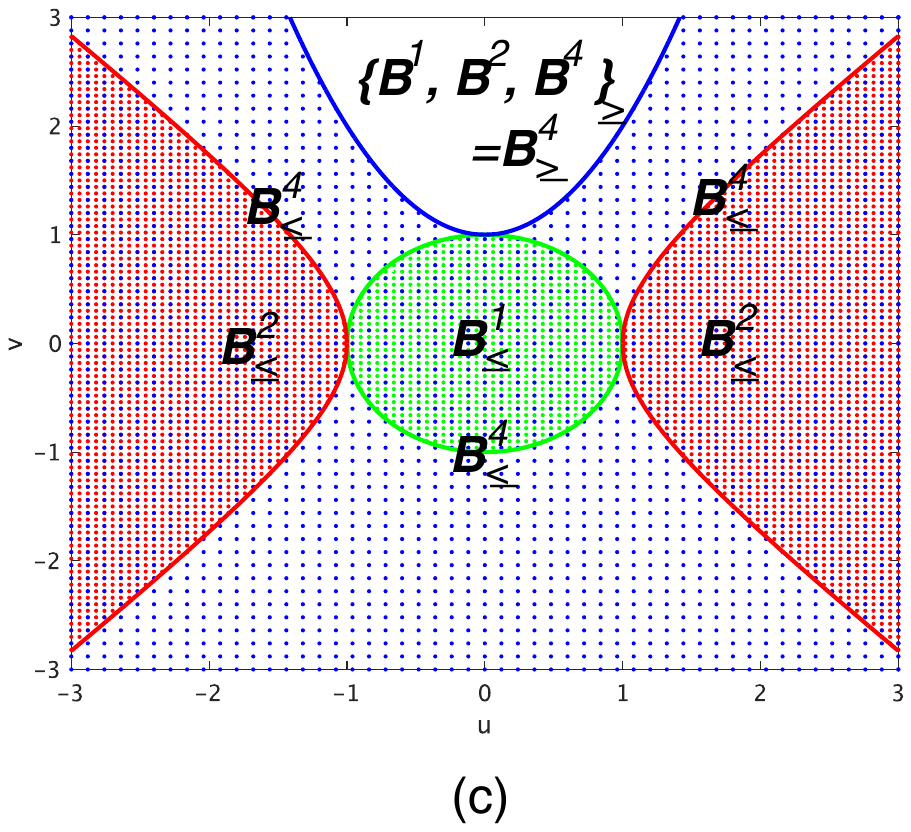}
\end{center}

\vspace{-20mm} 

\caption{
Illustration for conditions (C) and (D). 
$n=3$, $\coneF = \SymMat^3_+$ and $L_1(\coneF) = \Real^2$. 
(a) $\BC^1 = \{\B^1,\B^2,\B^3\}$, 
(b) $\BC^2 = \{\B^1,\B^3,\B^5\}$ and (c) $\BC^3 = \{\B^1,\B^2,\B^4\}$, where 
$q(\u,\B^1) = u_1^2 + u_2^2 - 1$,  
$q(\u,\B^2) = -u_1^2 + u_2^2 + 1$,  $q(\u,\B^3) = u_1^2-u_2 + 1$,  
$q(\u,\B^4) = -q(\u,\B^3)$, and $q(\u,\B^5) = u_2 + 1$. 
$\BC^k_\geq: $ the unshaded region $(k=1,2,3)$.  
} \label{fig:Fig1}
\end{figure} 

\medskip

Condition (B) (also (C) and (D)) depends on 
the choice of a face $\coneF$ of $\SymMat^n_+$ 
that  includes $\coneJ_+(\BC)$. If $\coneF^1$  and $\coneF^2$ are faces of 
$\SymMat^n_+$  such that $\coneF^1 \supseteq \coneF^2 \supseteq \coneJ_+(\BC)$,  
then the condition with $\coneF = \coneF^1$ 
implies the condition with $\coneF = \coneF^2$. 
Thus,  by choosing the minimal face 
$\coneF_{\min}$ of $\SymMat^n_+$ that contains $\coneJ_+(\BC)$, 
we obtain the weakest condition, 
since the condition with 
$\coneF_{\min}$ is implied by the condition with any larger face $\coneF$.
In this case, $\coneJ_+(\Phi^*(\BC))$ satisfies 
Slater's constraint qualification $\coneJ_+(\Phi^*(\BC)) \cap \SymMat^r_{++} \not= \emptyset$, 
where $\SymMat^r_{++}$ denotes the interior of $\SymMat^r_+$, {\it i.e.}, 
the set of $r \times r$ positive definite matrices.

\medskip

We now outline possible applications. 
QCQPs of the form~\eqref{eq:QCQPBCFH} arise in a variety of domains, 
including robotics and autonomous systems for avoiding multiple 
exclusion zones \cite{Hettich1993,Richards2002,Schulman2013},
sensor placement problems \cite{Boyd2004,Chong2003}
where optimal locations must lie outside risky regions, and data
classification tasks  involving the placement of test points outside known clusters \cite{Pang2021}.
The infeasible region of QCQP~\eqref{eq:QCQPBCFH}, given by the interior
of $\B_\leq$ $(\B \in \BC)$,  corresponds to the multiple exclusion zones,
risky regions, and known clusters, respectively 
(see Figures \ref{fig:Fig1}, \ref{fig:Fig3}, \ref{fig:Fig4} and \ref{fig:Fig5}). 
When the formulated QCQPs satisfy a NIQCC, 
such as conditions (I) and (III), they can be solved exactly through  their SDP relaxations. 
Furthermore, the proposed conditions (B), (C) and (D) are expected to substantially 
broaden the scope of such applications;  
for example,  by allowing the free addition of linear equality constraints and 
by effectively handling degenerate cases
 in which the feasible region fails to satisfy Slater’s constraint qualification.
Additionally, in Section~\ref{section:examples}, we illustrate 
some geometric examples 
for the reader interested in possible applications. 

\begin{figure}[t!]  
\begin{picture}(400,250)(0,0)
\put(145,245){I}
\put(133,240){\framebox(30,20)}
\put(220,245){II}
\put(210,240){\framebox(30,20)}
\put(330,245){III}
\put(320,240){\framebox(40,20)}
\put(70,220){generalization $\Downarrow$ }
\put(220,220){$\Downarrow$ Section~\ref{section:IItoB}}
\put(335,220){$\Downarrow$ generalization}
\put(17,195){(C)}
\put(10,190){\framebox(30,20)}
\put(70,195){$\Longrightarrow$}
\put(47,182){Theorem~\ref{theorem:main2}}
\put(170,195){(B)}
\put(120,190){\framebox(120,20)}
\put(270,195){$\Longleftarrow$}
\put(330,195){(D)}
\put(320,190){\framebox(40,20) }
\put(255,182){Section~\ref{section:DtoB}}
\end{picture}
\vspace{-60mm} 
\caption{Relationships among conditions (I), (II), (III), (B), (C) and (D). All conditions 
are equivalent under additional assumptions including 
Slater's constraint qualification 
($\coneJ_+(\BC)$ intersects with the interior of $\SymMat^n_+$ 
or $\coneF_{\min} = \SymMat^n_+$) 
and no redundancy on $\BC$ to represent $\coneJ_+(\BC)$. 
See Section~\ref{section:IandIIBC}. 
}
\end{figure}

\subsection{Outline of the paper}

In Section~\ref{section:exampleFBmain1},  a simple example is provided to demonstrate the   
effectiveness of the proposed conditions (B), (C) and (D) in comparison to 
the existing conditions (I), (II) and (III). In Section~\ref{section:facialReductionOfJ+}, we introduce 
an isomorphism $\Phi$ from a face $\coneF$ of $\SymMat^n_+$ onto $\SymMat^r_+$ 
and some basic theoretical issues related on $\Phi$ that play an essential role in the 
subsequent sections. Based on them, 
we  present a proof of Theorem~\ref{theorem:main1} 
in Section~\ref{section:ProofMain1}, and proofs of Theorem~\ref{theorem:main2}
and~\ref{theorem:JOYCE2024ext} in 
Section~\ref{section:Main2}. In Section~\ref{section:equivalence}, 
we investigate the relationships of the conditions presented above, (I), (II), (III), (B), (C), and (D) 
in detail. It is shown that 
condition (B) is the weakest among them, since it is implied by all the others,
and that all four conditions (I), (II), (B), and  (C) are equivalent under Slater’s constraint qualification 
and the absence of redundant constraints.  See Figure 2. 
Under these assumptions, condition (II) may be 
regarded as the dual of condition (B) (see Remark~\ref{remark:dualOfB}). 
In Section~\ref{section:examples}, we provide three geometric QCQP examples that 
satisfy condition (C).
Finally, Section~\ref{section:conclusion} contains our concluding remarks.

%\input sect2AnExample.tex
%!TEX root = ./main.tex

\section{An example illustrating the effectiveness of Theorems~\ref{theorem:main1}, 
\ref{theorem:main2} and~\ref{theorem:JOYCE2024ext} in 
comparison to Theorems~\ref{theorem:existing0} and~\ref{theorem:JOYCE2024}}

\label{section:exampleFBmain1}

 \begin{figure}[t!] 
\begin{center}
\includegraphics[height=60mm]{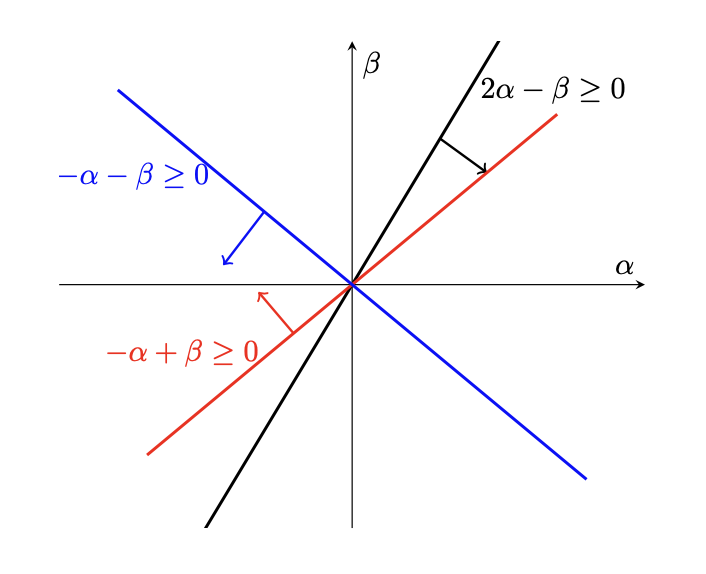}

\vspace{-5mm} 

\caption{The region of $(\alpha,\beta)$ determined by the inequalities
in~\eqref{eq:alphaBeta}.  
} 
\end{center}
\end{figure}

\medskip

\noindent
In this section, we compare our main results stated in Theorems~\ref{theorem:main1}, 
\ref{theorem:main2} and~\ref{theorem:JOYCE2024ext} 
with the known results  in Theorems~\ref{theorem:existing0} and~\ref{theorem:JOYCE2024} 
through an example to demonstrate the effectiveness of our main results. 
Let $n=4$, %$\H = \mbox{diag}(0,0,0,1)$, 
$\BC = \{\B^1,\B^2,\B^3\}$ and $\coneJ = \coneJ_+(\BC)$, where 
\begin{eqnarray*}
& & \B^1 = {\scriptsize \begin{pmatrix} 
-1 & -2 & 0 & 1 \\ 
-2 &  -1 & 0 & 0 \\
 0 & 0 & 2 & 2 \\ 
 1 & 0 & 2 & 2 
 \end{pmatrix}}, \ 
\B^2 = {\scriptsize \begin{pmatrix} 
1 &  1 & 0 & -2 \\ 
1 & -1 & 0 & 0  \\ 
0 & 0 & -1 & -1 \\ 
-2 & 0 & -1 & -1 \\ 
\end{pmatrix}}, \ 
\B^3 = {\scriptsize \begin{pmatrix}  
-1 & 1 & 0 & 1  \\ 
 1 & 1 & 0 & 0 \\ 
0 & 0 & -1 & -1 \\ 
1 & 0 & -1 & -1 
\end{pmatrix}}. 
\end{eqnarray*}
Let
\begin{eqnarray*}
\X^1 = {\scriptsize \begin{pmatrix}1&1&0&0\\1&1&0&0\\0&0&0&0\\0&0&0&0\end{pmatrix}}, \ 
\X^2 = {\scriptsize \begin{pmatrix} 1&0&0&1\\0&0&0&0\\0&0&0&0\\1&0&0&1\end{pmatrix}}, \ 
\X^3 = {\scriptsize \begin{pmatrix} 1&0&0&0\\0&0&0&0\\0&0&1&0\\0&0&0&0\end{pmatrix}}. 
\end{eqnarray*}
Then 
\begin{eqnarray*}
& & \X^1 \in \coneJ_+(\{\B^2,\B^3\})\backslash\coneJ_+(\{\B^1\}), \
\X^2 \in \coneJ_+(\{\B^1,\B^3\})\backslash\coneJ_+(\{\B^2\}), \\ 
& & \X^3 \in \coneJ_+(\{\B^1,\B^2\})\backslash\coneJ_+(\{\B^3\}).
\end{eqnarray*} 
This implies that none of $\B^1, \ \B^2$ and $\B^3$ is redundant to describe 
$\coneJ_+(\{\B^1,\B^2,\B^3\})$. 

\medskip

Let 
$\widehat{\X} = \mbox{diag} (1,1,0,0) $ \ (the $4 \times 4$ diagonal matrix with 
elements $1,1,0,0$). Then $\widehat{\X}  \in \coneJ_0(\B^2)$ 
but $\widehat{\X}  \not\in \coneJ_+(\B^1)$.
Hence $\BC$ does not satisfy 
condition (I) assumed in Theorem~\ref{theorem:existing0}. 
Also, condition (II) assumed in Theorem~\ref{theorem:existing0} dose not hold. In fact, assume on the contrary that 
(II) holds, which implies that the diagonal of  
$\alpha\B^1+\beta\B^2 \in \SymMat^4_+$ 
is nonnegative for some nonzero $(\alpha,\beta)\in \Real^2$. Hence 
\begin{eqnarray}
 -\alpha +\beta \geq 0, \ -\alpha-\beta \geq 0, \ 2\alpha-\beta \geq 0, \ 2\alpha-\beta \geq 0. 
\label{eq:alphaBeta}
\end{eqnarray} 
Clearly, only $(\alpha,\beta) = \0$ satisfies the above inequalities, 
as shown in  Figure 3, which 
 illustrates the region of $(\alpha,\beta)$ determined by the inequalities in~\eqref{eq:alphaBeta}. 
Therefore, condition (II) dose not hold. It is also easily verified that 
if $\overline{\u} = {\scriptsize \begin{pmatrix} 0 \\ \sqrt{2} \\ 0 \end{pmatrix}} \in \Real^3$ then $\bar{\u} \in \B^1_=$ but $\bar{\u} \not\in \B^2_\geq$. 
Therefore,~\eqref{eq:condOYCE2024} in condition (III) is not satisfied.

\medskip

To see whether condition (B) is satisfied for some face $\coneF$ that contains $\BC$, 
we observe that 
\begin{eqnarray*}
\coneJ_+(\BC) & \subseteq &  \left\{ \X\in \SymMat^4_+ : \inprod{\B^1+\B^2+\B^3}{\X} \geq 0 \right\} \\
& = & \left\{ \X\in \SymMat^4_+ : 
\Inprod{
{\scriptsize \begin{pmatrix} 
-1&0&0&0\\
0&-1&0&0\\
0&0&0&0\\
0&0&0&0
\end{pmatrix}}
}{\X} \geq 0 \right\}\\ 
& = & \left\{ 
{\scriptsize \begin{pmatrix}  \O & \O  \\ \O & \U \end{pmatrix}}
\in \SymMat^4: \U \in \SymMat^2_+ \right\}.
\end{eqnarray*}
Hence, by letting
\begin{eqnarray*}
& & \coneF^1 = \left\{ 
{\scriptsize \begin{pmatrix}  O & \O  \\ \O & \U \end{pmatrix}}
\in \SymMat^4: \U \in \SymMat^2_+ \right\}, \\  
& & \widetilde{\B^1} = {\scriptsize \begin{pmatrix} 2 & 2 \\ 2 & 2 \end{pmatrix}}, \ 
\widetilde{\B^2} = {\scriptsize \begin{pmatrix} -1 & -1  \\ -1 & -1 \end{pmatrix}}, \ 
\widetilde{\B^3} = {\scriptsize \begin{pmatrix}  -1 & -1 \\ -1 & -1 \end{pmatrix}}, \ 
\widetilde{\BC} = \{\widetilde{\B^1},\widetilde{\B^2},\widetilde{\B^3}\},  
\end{eqnarray*}
we see that 
$\coneJ_+(\BC) \subseteq  \coneF^1$, which indicates that $X_{ij} = 0$ $(1\leq i\leq 2 \ \mbox{or } 1\leq j\leq 2)$ if 
$\X \in \coneJ_+(\BC)$, 
and condition (B) is equivalent to 
\begin{eqnarray*}
\coneJ_0(\widetilde{\B}) \subseteq \coneJ_+(\widetilde{\A}) \ \mbox{or } 
\coneJ_+(\widetilde{\A}) \subseteq \coneJ_+(\widetilde{\B}) \ \mbox{for every } % distinct } 
\widetilde{\A}, \ \widetilde{\B} \in \widetilde{\BC}, 
\end{eqnarray*}
which  obviously holds since 
\begin{eqnarray*}
\coneJ_0(\widetilde{\B^1}) & = & \coneJ_0(\widetilde{\B^2}) = \coneJ_0(\widetilde{\B^3}) 
= \coneJ_+(\widetilde{\B^2}) = \coneJ_+(\widetilde{\B^3}) \\
 & = & \left\{ \U \in \SymMat^2_+: 
\inprod{{\scriptsize \begin{pmatrix} 1 & 1 \\ 1 & 1 \end{pmatrix}}}{\U} = 0\right\} 
\subseteq \SymMat^2_+ = \coneJ_+(\widetilde{\B^1}). 
\end{eqnarray*}

\medskip

Now, we verify whether condition (C) is satisfied with $\coneF = \coneF^1$. By definition, 
\begin{eqnarray*}
L(\coneF^1) & = & \left\{ \x \in \Real^4: \x\x^T \in \coneF^1 \right\} = \left\{ {\scriptsize \begin{pmatrix}0\\0\\x_3 \\x_4 \end{pmatrix}} \in \Real^4: x_3, \ x_4 \in \Real\right\}, \\
L_1(\coneF^1) & = & \left\{ \u \in \Real^3 : {\scriptsize \begin{pmatrix}\u\\1\end{pmatrix}}\in L(\coneF^1) \right\} 
= \left\{ {\scriptsize \begin{pmatrix}0\\0\\u_3\end{pmatrix}} \in \Real^3 : u_3 \in \Real \right\}. 
\end{eqnarray*}
For every $\u \in  L_1(\coneF^1)$, 
\begin{eqnarray*}
q(\u,\B^k) & = & 
\inprod{\widetilde{\B^k}}{{\scriptsize \begin{pmatrix}u_3\\1\end{pmatrix}}} = 
\left\{
\begin{array}{ll}
2(u_3+1)^2 & (k=1) \\
-(u_3+1)^2  & (k=2,3). 
\end{array}
\right. 
\end{eqnarray*}
holds. Hence 
\begin{eqnarray*}
& & \B^1_\leq \cap L_1(\coneF^1)  =  \B^2_\geq \cap L_1(\coneF^1) = \B^3_\geq \cap L(\coneF^1) = \left\{{\scriptsize \begin{pmatrix}0\\0\\-1\end{pmatrix}}\right\}, \\
& &  \B^1_\geq \cap L(\coneF^1) = \B^2_\leq \cap L_1(\coneF^1) = \B^3_\leq \cap L_1(\coneF^1) = 
\left\{{\scriptsize \begin{pmatrix}0\\0\\u_3\end{pmatrix}}: u_3 \in \Real\right\}. 
%, \\ 
%& & \B^2_\geq \cap L_1(\coneF^1) = \B^3_\geq \cap L_1(\coneF^1) = \left\{{\scriptsize \begin{pmatrix}0\\0\\-1\end{pmatrix}}\right\} \subseteq  \left\{{\scriptsize \begin{pmatrix}0\\0\\u_3\end{pmatrix}}: u_3 \in \Real\right\} = \B^1_\geq \cap L(\coneF^1).  
\end{eqnarray*}
%}
Therefore condition (C) with $\coneF = \coneF^1$ is satisfied. It is also easy to see %easily seen 
that condition (D) with $\coneF = \coneF^1$ is satisfied. 

\medskip

We note that 
$\coneF^1$ is not the minimal face of $\SymMat^4_+$ that contains $\coneJ_+(\BC)$. 
In fact, if we let 
\begin{eqnarray*}
\coneF^2 = \left\{ \X\in \SymMat^4_+ : 
\Inprod{{\scriptsize \begin{pmatrix}  0 & 0 & 0 & 0  \\ 0 & 0 & 0 & 0 \\ 0 & 0 & 1 & 1 \\ 0 & 0 & 1 & 1 \end{pmatrix}}}{\X} = 0 \right\},  
\end{eqnarray*}
then $\coneF = \coneF^1 \cap \coneF^2$ forms the minimal face.
In this case 
\begin{eqnarray*}
& & L(\coneF) = \left\{{\scriptsize \begin{pmatrix}0\\0\\x_3\\x_4\end{pmatrix}}:x_3+x_4=0\right\}, \ L_1(\coneF) = \left\{
{\scriptsize \begin{pmatrix}0\\0\\-1\end{pmatrix}}\right\},\\
& & \coneJ_+(\BC) \cap \coneF = \coneF = 
\left\{ {\scriptsize \begin{pmatrix}  O & \O  \\ \O & \U \end{pmatrix}}\in \SymMat^4: 
\inprod{{\scriptsize \begin{pmatrix}  1 & 1  \\ 1 & 1 \end{pmatrix}}}{\U} = 0, \ \U \in \SymMat^2_+ \right\}, \\
& & \BC_\geq \cap L_1(\coneF) = L_1(\coneF) = 
\left\{ {\scriptsize \begin{pmatrix}  0 \\ 0  \\ -1 \end{pmatrix}} \right\}
\end{eqnarray*}
hold.

%\input sect3Main1.tex
%!TEX root = ./main.tex
%\newpage

\section{On Theorem~\ref{theorem:main1}}

\label{section:Main1}

Throughout Sections~\ref{section:Main1} and~\ref{section:Main2}, 
let $\BC \subseteq \SymMat^n$ and $\coneF$ be a face of $\SymMat^n_+$ that 
contains $\coneJ_+(\BC)$. 

\subsection{
Facial reduction of $\coneJ_+(\BC)$ into $\SymMat^r_+$
}

\label{section:facialReductionOfJ+}

We first represent $\coneF$ as 
$\coneF = \left\{ \X \in \SymMat^n_+ : \inprod{\F}{\X} = 0 \right\}$ for some 
$\F \in \SymMat^n_+$. Suppose that $\mbox{rank}\F = n-r$ 
for some $r \in \{0,1,\ldots,n\}$. We can take 
an $n \times n$ orthogonal matrix $\P$ that diagonalizes $\F$ such that 
$\P^T\F\P = \mbox{diag}(0,\ldots,0,\lambda_{r+1},\ldots,\lambda_n) \in \SymMat^n_+$ 
for some positive numbers $\lambda_{r+1},\ldots,\lambda_n$, where $0,\ldots,0,\lambda_{r+1},\ldots,\lambda_n$ denote the eigenvalues of $\F$ with $\lambda_i > 0$ 
$(r+1 \leq i \leq n)$. Let $\P_{1\_r}$ denote the $n \times r$ matrix consisting of the first 
$r$ columns of $\P$; each $j$ th columns of $\P_{1\_r}$ is an 
eigenvector of $\F$ associated with  the zero eigenvalue. Let $\M$  be an arbitrary 
$r \times r$ nonsingular matrix.  Define 
\begin{eqnarray*}
\Phi(\X) & = & \M^T\P_{1\_r}^T\X\P_{1\_r}\M \in \SymMat^r_+ \ \mbox{for every } \X \in \coneF, \\
\Phi^*(\Y) & = & \M^{-1}\P_{1\_r}^T\Y\P_{1\_r}\M^{-T} \ \mbox{for every } \Y \in \SymMat^n, \\
\theta(\x) & = & \M^T\P_{1\_r}^T\x\ \mbox{for every } \x \in L(\coneF). 
\end{eqnarray*}
Then $\Phi: \coneF \rightarrow \SymMat^r_+$ forms a linear isomorphism from $\coneF$ onto 
$\SymMat^r_+$ and $\Phi^* : \SymMat^n \rightarrow \SymMat^r$ the adjoint map with respect to 
$\Phi$ (see, for example, \cite{BORWEIN1981,PATAKI2000,PATAKI2013}). Here 
$\V \in \SymMat^r_+ \rightarrow \M^T\V\M \in \SymMat^r_+$ serves as an 
automorphism on $\SymMat^r_+$. The choice of $\M$ is not relevant here, but it becomes 
relevant   in Section~\ref{section:Main2} where   
  QCQP~\eqref{eq:QCQPBCFH} is considered,  {\it i.e.}, 
COP($\coneG^n\cap\coneJ_+(\BC)\cap\coneF,\Q,\H$) with 
$\H = \mbox{diag}(0,\ldots,0,1) \in \SymMat^n_+$.
So we may assume that 
$\M$ is the $r \times r$ identity matrix in this section. 

%%%%%%%%%%
\lemm \label{lemma:isomorphismS} {\ }\vspace{-2mm}
\begin{description}
\item{(i) }
The map $\V \in \SymMat^r_+ 
\rightarrow \P_{1\_r}\M^{-T} \V\M^{-1}\P_{1\_r}^T \in \coneF$ serves as 
the inverse of $\Phi : \coneF \rightarrow \SymMat^r_+$. 
\vspace{-2mm} 
\item{(ii) }
$\theta: L(\coneF) \rightarrow \Real^r$ is linear, one-to-one and onto.\vspace{-2mm}
\item{(iii) } $\inprod{\A}{\X} = \inprod{\Phi^*(\A)}{\Phi(\X)}$ for every 
$\A \in \SymMat^n$ and $\X \in \coneF$.\vspace{-2mm}
\item{(iv) } $\inprod{\A}{\x\x^T} = \inprod{\Phi^*(\A)}{\theta(\x)\theta(\x)^T}$ for every 
$\A \in \SymMat^n$ and $\x \in L(\coneF)$.\vspace{-2mm}
\item{(v) } $\Phi(\coneJ_-(\B))\cap\coneF) = \coneJ_-(\Phi^*(\B)), \
\Phi(\coneJ_0(\B))\cap\coneF) = \coneJ_0(\Phi^*(\B))$ and 
$\Phi(\coneJ_+(\B))\cap\coneF) = \coneJ_+(\Phi^*(\B))$ for every 
$\B \in \BC$.\vspace{-2mm}
\item{(vi) } $\Phi(\coneJ_+(\BC))) = \Phi(\coneJ_+(\BC))\cap\coneF) = \coneJ_+(\Phi^*(\BC))$.
\end{description}
\elemm
%%%%%%%%%%
\proof{
(i) For every $\V \in \SymMat^r_+$,  we see 
\begin{eqnarray*}
\lefteqn{ \inprod{\F}{\P_{1\_r}\M^{-T} \V\M^{-1}\P_{1\_r}^T}} \\
%& = & 
%\inprod{\P\mbox{diag}(0,\ldots,0,\lambda_{r+1},\ldots,\lambda_n)\P^T}{
%\P_{1\_r}\M^{-T} \V\M^{-1}\P_{1\_r}^T} \\
& = &  \inprod{\P_{1\_r}^T\P\mbox{diag}(0,\ldots,0,\lambda_{r+1},\ldots,\lambda_n)\P^T\P_{1\_r}}{
\M^{-T} \V\M^{-1}} \\
& = & \inprod{\O}{\M^{-T} \V\M^{-1}} = 0. 
\end{eqnarray*}
%}
Hence $\P_{1\_r}\M^{-T} \V\M^{-1}\P_{1\_r}^T \in \coneF$. Also  
$
\Phi(\P_{1\_r}\M^{-T} \V\M^{-1}\P_{1\_r}^T) = \V
$ 
for every $\V \in \SymMat^r_+$. Therefore, the desired result follows. 
%
%\TBrown{
%\begin{eqnarray*}
%\Phi(\P_{1\_r}\M^{-T} \V\M^{-1}\P_{1\_r}^T) = \M^T\P_{1\_r}^T(\P_{1\_r}\M^{-T} \V\M^{-1}\P_{1\_r}^T)\P_{1\_r}\M = \V
%\end{eqnarray*}
%}
%
\medskip 

(ii) 
By definition, $\theta: L(\coneF) \rightarrow \Real^r$ is linear and 
$ \theta(L(\coneF)) \subseteq \Real^r$. Let $\v \in \Real^r$. Define 
$\x = \P_{1\_r}\M^{-T} \v$ and $\V = \v\v^T \in \SymMat^r_+$. 
Then $\x\x^T = \P_{1\_r}\M^{-T} \V \M^{-1}\P_{1\_r}^T \in \coneF$ 
as we have seen in the proof (i) above. Hence $\x \in L(\coneF)$ by 
the definition of $L(\coneF)$. 
We also see that $\theta(\x) = \M^T\P_{1\_r}^T\x = \M^T\P_{1\_r}^T\P_{1\_r}\M^{-T} \v = \v$. 
Hence $\theta(L(\coneF)) = \Real^r$, and we have shown $\theta(L(\coneF)) = \Real^r$. 
To see that $\theta: L(\coneF) \rightarrow \Real^r$ is one-to-one, assume that 
$\theta(\x^1) = \theta(\x^2)$ for some $\x^1, \x^2 \in L(\coneF)$. Then $\theta(\x^2-\x^1) = \0$.
Hence $\SymMat^r_+ \ni \O = \theta(\x^2-\x^1)\theta(\x^2-\x^1)^T = \Phi((\x^2-\x^1)(\x^2-\x^1)^T)$, 
which implies $(\x^2-\x^1)(\x^2-\x^1)^T = \O$ 
and $\x^2-\x^1 = \0$.

\medskip 

(iii) Let $\A \in \SymMat^n$, $\X \in \coneF$ and $\V = \Phi(\X) \in \SymMat^r_+$. By (i), $\X = \P_{1\_r}\M^{-T} \V\M^{-1}\P_{1\_r}^T \in \coneF$. 
Thus, it follows that 
\begin{eqnarray*}
& & \inprod{\A}{\X} = 
\inprod{\A}{\P_{1\_r}\M^{-T} \V\M^{-1}\P_{1\_r}^T} 
= \inprod{\M^{-1}\P_{1\_r}^T\A\P_{1\_r}\M^{-T} }{\V} = \inprod{\Phi^*(\A)}{\Phi(\X)}.
\end{eqnarray*}

\medskip

(iv) Let $\A \in \SymMat^n$ and 
$\x \in L(\coneF)$. Then $\x\x^T \in \coneF$. 
Consequently, $\inprod{\A}{\x\x^T} = \inprod{\Phi^*(\A)}{\Phi(\x\x^T)} 
= \inprod{\Phi^*(\A)}{\theta(\x)\theta(\x)^T}$ follows from (iii).

\medskip 

Assertion (v) can be proved easily by assertions (i) and (iii). Assertion (vi) 
follows from assertion (v). 
\qed
}

\medskip

The mapping $\coneJ_+(\BC)  \subseteq \coneF \rightarrow 
\Phi(\coneJ_+(\BC))\cap\coneF) = \coneJ_+(\Phi^*(\BC)) 
\subseteq \SymMat^r$ 
in (vi) is interpreted as {\em a facial reduction} \cite{BORWEIN1981}
 of $\coneJ_+(\BC)) \subseteq \coneF$ 
into $\SymMat^r$. 
For simplicity of notation, we denote $\Phi^*(\B)$ by $\widetilde{\B}$ for every $\B \in \SymMat^n$, 
and $\Phi^*(\BC) = \{\Phi^*(\B): \B \in \BC\}$ by $\widetilde{\BC}$. Then, 
we can simplify 
QCQP~\eqref{eq:QCQPBC0},  
its SDP relaxation~\eqref{eq:SDPBC},  
and condition (B) by the mappings 
$\Phi$ from $\coneF$ onto $\SymMat^r_+$ 
and $\theta$ from $L(\coneF)$
 onto $\Real^r$ as follows. 

\medskip

\noindent
COP($\coneG^r\cap\coneJ_+(\widetilde{\BC}),\widetilde{\Q},\widetilde{\H}$): 
\begin{eqnarray}
\eta(\coneG^r\cap\coneJ_+(\widetilde{\BC}),\widetilde{\Q},\widetilde{\H}) 
& = & \inf\left\{ \inprod{\widetilde{\Q}}{\v\v^T} : 
\begin{array}{l}
\v \in \Real^r, \
\v\v^T\in\coneJ_+(\widetilde{\BC}), \\
\inprod{\widetilde{\H}}{\v\v^T} = 1
\end{array}
\right\}.\label{eq:QCQP10}
\end{eqnarray}
COP($\coneJ_+(\widetilde{\BC}),\widetilde{\Q},\widetilde{\H}$): 
\begin{eqnarray}
\eta(\coneJ_+(\widetilde{\BC}),\widetilde{\Q},\widetilde{\H}) 
& = & \inf\left\{ \inprod{\widetilde{\Q}}{\V} : 
\begin{array}{l}
\V \in \SymMat^r_+, \
\V\in\coneJ_+(\widetilde{\BC}), \\
\inprod{\widetilde{\H}}{\V} = 1
\end{array}
\right\}.\label{eq:SDP10}
\end{eqnarray}
\begin{description}
\item{($\widetilde{\rm \B}$) } For every % distinct 
$\widetilde{\A}, \widetilde{\B} \in \widetilde{\BC}$, 
either $\coneJ_0(\widetilde{\B}) \subseteq \coneJ_+(\widetilde{\A})$ 
or $\coneJ_+(\widetilde{\A}) \subseteq \coneJ_+(\widetilde{\B})$ holds. 
\end{description}

\medskip

\noindent
More precisely, Lemma~\ref{lemma:isomorphismS} ensures: 
\begin{description}
\item{[a] } 
$\x \in \Real^n$ is a feasible solution of~\eqref{eq:QCQPBC0} with the objective value 
$\inprod{\Q}{\x\x^T}$ if and only if $\v = \theta(\x)$ 
is a feasible solution of~\eqref{eq:QCQP10} with the objective value 
$\inprod{\widetilde{\Q}}{\v\v^T} = \inprod{\Q}{\x\x^T}$.
\item{[b] } 
$\X \in \SymMat^n_+$ is a feasible solution of~\eqref{eq:SDPBC} with the objective value 
$\inprod{\Q}{\X}$ if and only if $\V = \Phi(\X)$ 
is a feasible solution of~\eqref{eq:SDP10} with the objective value 
$\inprod{\widetilde{\Q}}{\V} = \inprod{\Q}{\X}$.\vspace{-2mm}  
\item{[c] } Condition (B) holds if and only if condition ($\widetilde{\rm B}$) does.\vspace{-2mm} 
\end{description}
In addition, we know that when $\X \in \coneF$ and $\V = \Phi(\X) = \M^T\P_{1\_r}^T\X\P_{1\_r}\M$, 
which imply $\X = \P_{1\_r}\M^{-T} \V\M^{-1}\P_{1\_r}^T$, $\X$ is rank-1 if and only if so is $\V$. 
It follows that \vspace{-2mm}
\begin{description}
\item{[d] } $\coneJ_+(\BC) = \coneJ_+(\BC) \cap \coneF \in \wFC(\coneG^n)$ if and only if 
$\coneJ_+(\widetilde{\BC}) = \coneJ_+(\Phi^*(\BC)) = 
\Phi(\coneJ_+(\BC) \cap \coneF) 
\in \wFC(\coneG^r)$.  
\end{description}

\subsection{Proof of Theorem~\ref{theorem:main1}}

\label{section:ProofMain1}

In view of the discussion of Section~\ref{section:facialReductionOfJ+}, 
it suffices to prove $\coneJ_+(\widetilde{\BC}) \in \wFC(\coneG^r)$ 
under condition ($\widetilde{\rm B}$) (see [c] and [d]). 
For simplicity of notation, we omit  $\widetilde{\ }$ \ from $\widetilde{\B}$ 
$(\B \in \BC)$, $\widetilde{\BC}$, $\widetilde{\Q}$ and $\widetilde{\H}$ 
for the proof, or equivalently, we prove $\coneJ_+(\BC) \in \wFC(\coneG^n)$ 
under condition (B) with $\coneF = \SymMat^n_+$. 

\label{section:ProofMain1}

\subsubsection{The finite case of $\BC$ }

\label{section:ProofMain1Finite}

This case can  be derived easily from Theorem~\ref{theorem:existing0}. 
Assume that $\BC$ is finite and condition (B) holds with $\coneF = \SymMat^n_+$. 
If $\coneJ_+(\A)\subseteq  \coneJ_+(\B)$ for distinct 
$\A, \ \B \in \BC$, then $\B$ is redundant to 
describe $\coneJ_+(\BC)$, {\it i.e.}, $\coneJ_+(\BC\backslash \{\B\}) = \coneJ_+(\BC)$. 
Hence, we can remove such elements $\B$ from $\BC$ one by one, recursively, % from 
%$\BC$ recursively 
to construct a reduced set $\BC$. 
Then, the resulting $\BC$ eventually satisfies condition (I).
\qed

\subsubsection{The infinite case of $\BC$}

We present three lemmas for the proof. 

%%%%%%%%%%
\lemm \label{lemma:sequence}  
Let 
$\{\BC^k \subseteq \SymMat^n : k=1,2,\ldots \}$ 
be a sequence 
such that $\bigcap_{k=1}^\infty \coneJ_+(\BC^k) = \coneJ$ for some closed convex cone 
$\coneJ \subseteq \SymMat^n_+$. Then 
$\bigcap_{m=1}^\infty \mbox{co}\big(\coneG^n\cap\big(\bigcap_{k=1}^m \coneJ_+(\BC^k)\big)\big) 
= \mbox{co}\big(\coneG^n\cap\coneJ \big) $.  \vspace{-2mm}
\elemm
%%%%%%%%%%
\proof{Since $\bigcap_{k=1}^m \coneJ_+(\BC^k) \supseteq \coneJ$ $(m=1,2,\ldots)$, 
%
%\TBrown{
%$\mbox{co}\big(\coneG^n\cap\big(\bigcap_{k=1}^m \coneJ_+(\BC^k)\big)\big) 
%\supseteq \mbox{co}\big(\coneG^n\cap\coneJ \big)$ $(m=1,2,\ldots)$ and 
%}
%
$\bigcap_{m=1}^\infty \mbox{co}\big(\coneG^n\cap\big(\bigcap_{k=1}^m \coneJ_+(\BC^k)\big)\big) 
\supseteq \mbox{co}\big(\coneG^n\cap\coneJ \big)$ follows. To prove the converse inclusion, 
let $\overline{\X} \in \mbox{co}\big(\coneG^n \cap \big(\bigcap_{k=1}^{m} \coneJ_+(\BC^k)\big)\big)$ $(m=1,2,\ldots \ )$.  Then,  for each $m=1,2,\ldots \ $, 
there exist $\X_m^p \in\coneG^n \cap \big(\bigcap_{k=1}^{m} \coneJ_+(\BC^k)\big) \subseteq \SymMat^n_+$   
$(p=1,2,\ldots,\ell)$ for some $\ell \leq \mbox{dim}\SymMat^n = 
n(n-1)/2$ such that 
$\overline{\X} = \sum_{p=1}^{\ell}  \X_m^p$.
Let $q \in \{1,2,\ldots,\ell\}$ be fixed arbitrarily.  
Consider the sequence $\{ \X_m^q \in \SymMat^n_+: m=1,2,\ldots, \ \}$.  
The sequence is bounded since $\overline{\X} \in \SymMat^n_+$, 
$\X_m^p \in \SymMat^n_+$ $(p=1,\ldots,\ell)$ and 
$\inprod{\I}{\overline{\X}} = \inprod{\I}{\sum_{p=1}^{\ell}  \X_m^p} \geq  \inprod{\I}{\X_m^q}$. 
Hence the sequence admits a subsequence converging to some $\overline{\X}^q \in \SymMat^n_+$. 
For notational simplicity, we relabel this subsequence as the sequence itself. 
Then, 
$\overline{\X}^q \in \mbox{cl}\big(\bigcap_{k=1}^\infty \coneJ_+(\BC^k)) = \coneJ$.  
Since 
$\coneG^n$ is closed, we also see $\overline{\X}^q \in \coneG^n$. Therefore, taking the limit 
of the identity $\overline{\X} = \sum_{p=1}^{\ell}  \X_m^p$ as $m \rightarrow \infty$, we obtain that
$ 
\overline{\X}  = \sum_{p=1}^{\ell}  \overline{\X}^p \ \mbox{and } \overline{\X}^p \in \coneG^n \cap \coneJ \ (p=1,\ldots,\ell). 
$ 
Therefore, we have shown that 
$\overline{\X} \in \mbox{co}\big(\coneG^n\cap\coneJ\big)$ and 
$\bigcap_{m=1}^\infty \mbox{co}\big(\coneG^n\cap\big(\bigcap_{k=1}^m \coneJ_+(\BC^k)\big)\big) 
\subseteq \mbox{co}\big(\coneG^n\cap\coneJ \big)$. 
\qed
}

\smallskip

We may assume without loss of generality that $\BC$ is bounded since we can 
replace $\BC$ with $\BC' = \{\B/\parallel \B \parallel : \O \not=\B \in \BC\}$ if $\BC$ is unbounded, 
where $\parallel \B \parallel $ denotes the Frobenius norm of $\B \in \SymMat^n$.
For each $\epsilon > 0$, define an open neighborhood 
$U(\B,\epsilon) = \{\A \in \SymMat^n : \parallel \A - \B \parallel < \epsilon \}$ of each 
$\B \in \mbox{cl}\BC$. 
Let $\{\epsilon_k\}$ be a sequence of positive numbers that converges to $0$. Let $k$ be fixed. 
Since $\BC$ is bounded, cl$\BC$ is a compact in $\SymMat^n$. 
Hence we can take a finite subset $\BC^k \subseteq \BC$ such that the union of all $U(\B,\epsilon_k)$ $(\B \in \BC^k)$ covers $\BC$. 
Then, the sequence $\{\BC^k \ (k=1,2,\ldots) \}$ satisfies that 
\begin{eqnarray}
& & \BC^{k} \subseteq \BC,  \  \coneJ_+(\BC^k) \supseteq \coneJ_+(\BC) \ (k=1,2,\ldots), \ 
\bigcap_{k=1}^\infty \coneJ_+(\BC^k) \supseteq \coneJ_+(\BC),  
\label{eq:subset} \\
& & \forall \B \in \BC, \ \exists \B' \in \BC^{k}; \parallel \B' - \B \parallel < \epsilon_k \ (k=1,2,\ldots). 
\label{eq:distance} 
\end{eqnarray}

%%%%%%%%%%
\lemm \label{lemma:coneJ} 
$\bigcap_{k=1}^\infty \coneJ_+(\BC^k) = \coneJ_+(\BC)$.
\vspace{-2mm}
\elemm
%%%%%%%%%%
\proof{By~\eqref{eq:subset}, it suffices  to show that 
$\bigcap_{k=1}^\infty \coneJ_+(\BC^k) \subseteq \coneJ_+(\BC)$.  Let 
$\overline{\X} \in \bigcap_{k=1}^\infty \coneJ_+(\BC^k)$. To prove $\overline{\X} \in \coneJ_+(\BC)$, 
we show that $\inprod{\overline{\B}}{\overline{\X}} \geq 0$ 
for an arbitrarily chosen $\overline{\B} \in \BC$. 
By~\eqref{eq:distance}, there exists a sequence $\{\B_k \in \BC^k\}$ that converges 
$\overline{\B} \in \BC$ as $k \rightarrow \infty$. Since $\overline{\X} \in \coneJ_+(\BC^k)$, we see that
$ 
\inprod{\B_k}{\overline{\X}} \geq 0 \ (k=1,2,\ldots \ ) 
$. 
Hence, we obtain $\inprod{\overline{\B}}{\overline{\X}} \geq 0$ by taking the limit 
as $k \rightarrow \infty$. 
\qed
}

%%%%%%%%%%
\lemm \label{lemma:finite}
Assume that $\BC$ satisfies condition (B). Then 
$\mbox{co}\big(\coneG^n\cap\big(\bigcap_{k=1}^m \coneJ_+(\BC^k)\big)) = \bigcap_{k=1}^m \coneJ_+(\BC^k)$  $(m=1,2,\ldots )$.\vspace{-2mm}
\elemm 
%%%%%%%%%%
\proof{Let $m \in \{1,2,\ldots \}$ be fixed. 
We note that $\bigcap_{k=1}^m \coneJ_+(\BC^k) = \coneJ_+(\bigcup_{k=1}^m \BC^k)$. 
Each $\bigcup_{k=1}^m \BC^k$ satisfies condition (B) since it is a subset of 
$\BC$. Since $\bigcup_{k=1}^m \BC^k$ is finite, 
 $\coneJ_+(\bigcup_{k=1}^m \BC^k) \in \wFC(\coneG^n)$ 
 as shown in Section~\ref{section:ProofMain1Finite}. 
 Therefore, 
$\mbox{co}\big(\coneG^n\cap\big(\bigcap_{k=1}^m \coneJ_+(\BC^k)\big)) 
= \bigcap_{k=1}^m \coneJ_+(\BC^k)$.  
\qed
}

\medskip

Now, we show 
$\mbox{co}\big(\coneG^n \cap \coneJ_+(\BC)\big) = \coneJ_+(\BC)$, which is 
equivalent to $ \coneJ_+(\BC) \in \wFC(\coneG^n)$ by definition. 
By Lemmas~\ref{lemma:coneJ} and ~\ref{lemma:finite}, we see that
\begin{eqnarray*}
\mbox{co}\big(\coneG^n\cap\big(\bigcap_{k=1}^m \coneJ_+(\BC^k)\big))  
= \bigcap_{k=1}^m \coneJ_+(\BC^k) \supseteq \bigcap_{k=1}^\infty \coneJ_+(\BC^k) = 
\coneJ_+(\BC) \ 
(m=1,2,\ldots), 
\end{eqnarray*}
which implies that 
 $\bigcap_{m=1}^\infty \mbox{co}\big(\coneG^n\cap\big(\bigcap_{k=1}^m \coneJ_+(\BC^k)\big)) 
\supseteq \coneJ_+(\BC)$.  By Lemma~\ref{lemma:sequence}, 
$\bigcap_{m=1}^\infty \mbox{co}\big(\coneG^n\cap\big(\bigcap_{k=1}^m \coneJ_+(\BC^k)\big)\big)  = 
\mbox{co}\big(\coneG^n\cap\coneJ_+(\BC)\big)$. 
Therefore, we have shown that $\mbox{co}\big(\coneG^n\cap\coneJ_+(\BC)\big) \supseteq \coneJ_+(\BC)$. 
The converse inclusion 
$\mbox{co}\big(\coneG^n \cap \coneJ_+(\BC)\big) \subseteq \coneJ_+(\BC)$ is straightforward 
since $\coneJ_+(\BC)$ is convex. 
\qed

%
%\input sect4Main2.tex
%!TEX root = ./main.tex

\section{On Theorems~\ref{theorem:main2} and~\ref{theorem:JOYCE2024ext} }

\label{section:Main2}

\subsection{Facial reduction of $\BC_\geq\cap\L_1(\coneF)$ into $\Real^{r-1}$}

\label{section:facialFBmain2}

To adapt the argument in Section~\ref{section:facialReductionOfJ+} to 
QCQP~\eqref{eq:QCQPBCFH} ({\it i.e.}, 
COP($\coneG^n\cap\coneJ_+(\BC),\Q,\H^1$) with 
$\H^1 = \mbox{diag}(0,\ldots,0,1)\in \SymMat^n_+$), its SDP relaxation, 
conditions (C) and (D), we need some additional arguments.  First, we assume that the feasible 
region $\BC_\geq \cap L_1(\coneF)$ of 
QCQP~\eqref{eq:QCQPBCFH} 
is nonempty, as required  in condition (C). Consequently, the feasible region of 
COP($\coneG^r\cap\coneJ_+(\widetilde{\BC}),\widetilde{\Q},\widetilde{\H^1}$), 
QCQP \eqref{eq:QCQP10} with $\widetilde{\H} = \widetilde{\H^1}$
% } 
is also nonempty.
For any choice of a nonsingular $\M$, the rank of 
$\widetilde{\H^1} = \Phi^*(\H^1) = \M^{-1}\P_{1.r}^T\H^1\P_{1.r}\M^{-T} \in \SymMat^r_+$ is $1$, otherwise 
the rank is $0$ or $\widetilde{\H^1} = \O \in \SymMat^r_+$; hence  
COP($\coneG^r\cap\coneJ_+(\widetilde{\BC}),\widetilde{\Q},\widetilde{\H^1}$) 
\eqref{eq:QCQP10} is 
infeasible. Specifically, we can take a nonsingular matrix $\M$ such that 
$\widetilde{\H^1} = \M^{-1}\P_{1.r}^T\H^1\P_{1.r}\M^{-T} = \mbox{diag}(0,\ldots,0,1) \in 
\SymMat^r_+$. 
In this case, if $\u \in L_1(\coneF)$, then ${\scriptsize \begin{pmatrix} \u \\ 1 \end{pmatrix}} \in L (\coneF)$ and 
$
\inprod{\widetilde{\H^1}}{\theta \Big({\scriptsize \begin{pmatrix} \u \\ 1 \end{pmatrix}} \Big)
\theta\Big({\scriptsize \begin{pmatrix} \u \\ 1 \end{pmatrix}}\Big)^T}  =\inprod{\H^1}{{\scriptsize \begin{pmatrix} \u \\ 1 \end{pmatrix}}{\scriptsize \begin{pmatrix} \u \\ 1 \end{pmatrix}}^T} = 1 
$
hold by Lemma~\ref{lemma:isomorphismS} (iv). 
Since $\theta : L(\coneF) \rightarrow \Real^{r-1}$ is linear, we have either 
$\theta\Big({\scriptsize \begin{pmatrix} \u \\ 1 \end{pmatrix}}\Big)_r = 1$ for every $\u \in L(\coneF)$ or 
$\theta\Big({\scriptsize \begin{pmatrix} \u \\ 1 \end{pmatrix}}\Big)_r = -1$ for every $\u \in L(\coneF)$. 
Define $\theta_1 : L_1(\coneF) \rightarrow \Real^{r-1}$ by 
\begin{eqnarray*}
\theta_1(\u) = 
\left\{
\begin{array}{ll}
\mbox{
\Big( $\theta\Big({\scriptsize \begin{pmatrix} \u \\ 1 \end{pmatrix}}\Big)_1,\ldots,
\theta\Big({\scriptsize \begin{pmatrix} \u \\ 1 \end{pmatrix}}\Big)_{r-1}$}  \Big)^T
& \mbox{in the former case,} \\ 
\mbox{
$\Big( -\theta\Big({\scriptsize \begin{pmatrix} \u \\ 1 \end{pmatrix}}\Big)_1,\ldots,
-\theta\Big({\scriptsize \begin{pmatrix} \u \\ 1 \end{pmatrix}}\Big)_{r-1}$}  \Big)^T
& 
\mbox{in the latter case.} 
\end{array}
\right.
\end{eqnarray*}
Then, 
\begin{eqnarray*}
q(\u,\B) = \inprod{\B}{
{\scriptsize \begin{pmatrix} \u \\ 1 \end{pmatrix}}
{{\scriptsize \begin{pmatrix} \u \\ 1 \end{pmatrix}}^T} } = \inprod{\widetilde{\B}}{\theta\Big({\scriptsize \begin{pmatrix} \u \\ 1 \end{pmatrix}}\Big)\theta\Big({\scriptsize \begin{pmatrix} \u \\ 1 \end{pmatrix}}\Big)^T}
=  \inprod{\widetilde{\B}}{{\scriptsize \begin{pmatrix} \theta_1(\u) \\ 1 \end{pmatrix}}{\scriptsize \begin{pmatrix} \theta_1(\u) \\ 1 \end{pmatrix}}^T}
\end{eqnarray*}
for every $\u \in L_1(\coneF)$ and $\B \in \SymMat^n$. 
Now, for every 
$\w \in \Real^{r-1}$ 
and $\widetilde{\B} \in \SymMat^r$, 
we define 
\begin{eqnarray*}
\begin{array}{l}
q(\w,\widetilde{\B}) = 
{\scriptsize \begin{pmatrix} \w \\ 1 \end{pmatrix}}^T\widetilde{\B}
{\scriptsize \begin{pmatrix} \w \\ 1 \end{pmatrix}}, \\[3pt]
\widetilde{\B}_\geq, \ \widetilde{\B}_= \ \mbox{or } \widetilde{\B}_\leq  =  
\left\{ \w \in \Real^{r-1} :  q(\w,\widetilde{\B}) \geq 0,  \ = \ \mbox{or } \leq 0, \ \mbox{resp.}
\right\}, \\[3pt]
\widetilde{\BC}_\geq = \bigcap_{\widetilde{\B} \in \widetilde{\BC}} \widetilde{\B}_\geq = 
\left\{ \w \in \Real^{r-1} :  q(\w,\widetilde{\B}) \geq 0 \ (\widetilde{\B} \in \widetilde{\BC})\right\}.
\end{array}
\end{eqnarray*}
Using the notation above, QCQP~\eqref{eq:QCQPBCFH} is transformed into 
COP($\coneG^r\cap\coneJ_+(\widetilde{\BC}),\widetilde{\Q},\widetilde{\H^1}$) with 
$\widetilde{\H^1} = \mbox{diag}(0,\ldots,0,1) \in \SymMat^r_+$: 
\begin{eqnarray}
\eta(\coneG^r\cap\coneJ_+(\widetilde{\BC}),\widetilde{\Q},\widetilde{\H^1})
& = &  \inf\{q(\w,\widetilde{\Q}) : q(\w,\widetilde{\B}) \geq 0 \ (\widetilde{\B} \in \widetilde{\BC}) \} \nonumber \\ 
& = &  \inf\{q(\w,\widetilde{\Q}) : \w \in \widetilde{\BC}_\geq \}, 
\label{eq:QCQP11}
\end{eqnarray}
and conditions (C) and (D) into 
\begin{description}
\item{($\widetilde{\rm C}$)} $\widetilde{\BC}_\geq\not= \emptyset$, and 
$\emptyset \not= \widetilde{\B}_\leq \subseteq \widetilde{\A}_\geq$ or 
$\widetilde{\A}_\geq \subseteq \widetilde{\B}_\geq$
holds for every 
$\widetilde{\A}, \ \widetilde{\B}  \in \widetilde{\BC}$. \vspace{-2mm}
\item{($\widetilde{\rm D}$) } 
$\BC$ is finite.  
$q(\cdot,\widetilde{\B}) : \Real^{r-1} \rightarrow \Real$ $(\widetilde{\B} \in \widetilde{\BC})$ is not affine 
and 
$ 
\widetilde{\B}_=  \subseteq \widetilde{\A}_\geq \ \mbox{ for every  } \widetilde{\A}, \ \widetilde{\B} \in \widetilde{\BC}. 
$
%\vspace{-2mm}
\end{description}

%%%%%%%%%
\lemm \label{lemma:isomorphismQ} \mbox{ \ } \vspace{-2mm} 
\begin{description}
\item{(i) } $\theta_1 : L_1(\coneF) \rightarrow \Real^{r-1}$ is one-to-one and onto.\vspace{-2mm}
\item{(ii) } $q(\u,\A) = q(\theta_1(\u),\widetilde{\A}) \ \mbox{for every } \u \in L_1(\coneF) \ 
\mbox{and } \A \in \SymMat^n$. \vspace{-2mm}
\item{(iii) } $\theta_1(\B_\geq \cap L_1(\coneF)) = \widetilde{\B}_\geq, 
\ \theta_1(\B_= \cap L_1(\coneF)) = \widetilde{\B}_= \ 
\mbox{and } \theta_1(\B_\leq \cap L_1(\coneF)) = \widetilde{\B}_\leq$. 
\end{description}
\elemm
%%%%%%%%%
\proof{
(i) This assertion follows from Lemma~\ref{lemma:isomorphismS} (ii) and the definition of 
$\theta_1 : L_1(\coneF) \rightarrow \Real^{r-1}$.

\medskip

(ii) Let $\u \in L_1(\coneF)$ and   $\A \in \SymMat^n$ 
Then, 
\begin{eqnarray*}
q(\u,\A) & = & \inprod{\A}{{\scriptsize \begin{pmatrix}\u\\1\end{pmatrix}}
{\scriptsize \begin{pmatrix}\u\\1\end{pmatrix}}^T} \\
& = & \inprod{\widetilde{\A}}
{\theta\Big({\scriptsize \begin{pmatrix}\u\\1\end{pmatrix}}\Big)\theta\Big({\scriptsize \begin{pmatrix}\u\\1\end{pmatrix}}\Big)^T} \ \mbox{(by Theorem~\ref{lemma:isomorphismS} (iv))
} \\ 
& = & \inprod{\widetilde{\A}}
{{\scriptsize \begin{pmatrix}\theta_1(\u)\\1\end{pmatrix}}
{\scriptsize \begin{pmatrix}\theta_1(\u)\\1\end{pmatrix}}^T} 
\ = \ q(\widetilde{\A},\theta_1(\u)). 
\end{eqnarray*}

\medskip

(iii) By (i) and (ii), we see that
\begin{eqnarray*}
\theta_1(\B_\geq \cap L_1(\coneF)) & = & \left\{
\theta_1(\u) : q(\u,\B) \geq 0, \ \u \in L_1(\coneF) \right\} \\
& = & \left\{
\theta_1(\u) : q(\theta_1(\u),\widetilde{\B}) \geq 0, \ \theta_1(\u) \in \theta_1(L_1(\coneF)) \right\} \\
& = & \left\{
\w \in \Real^{r-1}: q(\w,\widetilde{\B}) \geq 0 \right\} 
= \widetilde{\B}_\geq.  
\end{eqnarray*}
The second and third identity can be proved similarly. 
\qed
}

\medskip

As a result of Lemma~\ref{lemma:isomorphismQ}, we obtain that: 
\vspace{-1mm} 
\begin{description}
\item{[e] } 
$\u\in \Real^{n-1}$ is a feasible solution of QCQP~\eqref{eq:QCQPBCFH} with the 
objective value $q(\u,\Q)$ if and only if $\w = \theta_1(\u)$ 
is a feasible solution of QCQP~\eqref{eq:QCQP11} with the objective value 
$q(\w,\widetilde{\Q}) = q(\u,\Q)$.\vspace{-2mm}
\item{[f]}  Condition (C) holds if and only if condition ($\widetilde{\rm C}$) holds.\vspace{-2mm}
\item{[g]}  Condition (D) holds if and only if condition ($\widetilde{\rm D}$) holds. 
\end{description}

\subsection{Proof of Theorem~\ref{theorem:main2}}

\label{section:ProofMain2}

We have shown in [c] of Section~\ref{section:facialReductionOfJ+} and [f] above 
that conditions (B) and (C) 
are equivalent to conditions ($\widetilde{\rm B}$)  and ($\widetilde{\rm C}$), respectively. 
Hence, for Theorem~\ref{theorem:main2}, it is sufficient to prove that 
condition ($\widetilde{\rm C}$)  implies condition ($\widetilde{\rm B}$). 
For simplicity of notation, we omit $\widetilde{\ }$ from $\widetilde{\B}$ 
$(\B \in \BC)$, $\widetilde{\BC}$, $\widetilde{\Q}$ and $\widetilde{\H^1}$,  
or equivalently, we prove that condition (C)  implies condition (B)
under the assumption that $\coneF = \SymMat^n_+$, $L(\coneF) = \Real^r$ and 
$L_1(\coneF)=\Real^{r-1}$. 

\medskip

%%%%%%%%%
\lemm \label{lemma:NScond00} 
 Let $\A, \ \B \in \SymMat^n$. Then 
\vspace{-2mm}
\begin{description}
\item{(i) } $\coneJ_-(\B) \subseteq \coneJ_+(\A) \Rightarrow  \coneJ_0(\B) \subseteq \coneJ_+(\A) $.
\vspace{-1mm}
\item{(ii) } \vspace{-9mm}
\begin{eqnarray}
& \hspace{-71mm} 
\coneJ_-(\B) \subseteq \coneJ_+(\A) 
\Leftrightarrow 
\coneG^n\cap\coneJ_-(\B) \subseteq \coneG^n\cap\coneJ_+(\A)
& \nonumber \\ 
&\hspace{-36mm}  \Updownarrow &  \nonumber  \\
& \hspace{-36mm} 
\{\x \in \Real^n : \inprod{\B}{\x\x^T} \leq 0 \} \subseteq 
\{\x \in \Real^n : \inprod{\A}{\x\x^T} \geq 0\} 
&    
\label{eq:AKK01}\vspace{-2mm}
\end{eqnarray}
\end{description}
\elemm 
%%%%%%%%%
\proof{
(i) is obvious since $\coneJ_0(\B) \subseteq \coneJ_-(\B)$. (ii) 
For the first $\Leftrightarrow$, $\Rightarrow$ is straightforward.  
Conversely, if $\coneG^n\cap\coneJ_-(\B) \subseteq \coneG^n\cap\coneJ_+(\A)$, then 
\begin{eqnarray*}
\coneJ_-(\B) = \mbox{co}(\coneG^n\cap\coneJ_-(\B)) \subseteq \mbox{co}(\coneG^n\cap\coneJ_+(\A)) = \coneJ_+(\A) 
\end{eqnarray*}
follows by Theorem~\ref{theorem:existing0}. 
Hence $\Rightarrow$ follows. The second $\Updownarrow$ 
is straightforward. 
\qed
}

\bigskip

%%%%%%%%%
\lemm 
\label{lemma:NScond03} 
Let $\A, \ \B \in \SymMat^n$.  
Assume that $\emptyset \not=  \B_\leq \subseteq \A_\geq$. 
Then $\coneJ_-(\B) \subseteq  \coneJ_+(\A)$. 
\vspace{-2mm}
\elemm
%%%%%%%%%
\proof{
By Lemma~\ref{lemma:NScond00} (ii), it suffices to prove~\eqref{eq:AKK01} 
under the assumption 
$\emptyset \not=  \B_\leq \subseteq \A_\geq$.  
For every subset $D\subseteq \Real^{n-1}$, define
\[
\coneH(D)=\left\{\binom{u}{\xi}\in\Real^n:\ u\in\Real^{n-1},\ \xi>0,\ u/\xi\in D\right\}.
\]
Since $\B_{\leq}\neq\emptyset$ and $\A_{\geq}\neq\emptyset$, we may apply the homogenization identity in \cite[Lemma~2]{STURM2003} to obtain
\begin{eqnarray*}
& & \{\x\in\Real^n:\langle \B,\x\x^T\rangle\le 0\}=\coneH(\B_{\le})\cup \coneH(-\B_{\le}),\\
& & \{\x\in\Real^n:\langle \A,\x\x^T\rangle\ge 0\}=\coneH(\A_{\ge})\cup \coneH(-\A_{\ge}).
\end{eqnarray*}
Since $ \B_{\le}\subseteq \A_{\ge}$, we obtain
$\coneH(\B_{\le})\subseteq \coneH(\A_{\ge})$ and $\coneH(-\B_{\le})\subseteq \coneH(-\A_{\ge})$.
Therefore,
\[
\coneH(\B_{\le})\cup \coneH(-\B_{\le}) \subseteq \coneH(\A_{\ge})\cup \coneH(-\A_{\ge}),
\]
which implies (16). 
\qed 
}

%%%%%%%%%
\rema The following example shows the necessity of the assumption 
$\emptyset \not= \B_\leq$ in Lemma~\ref{lemma:NScond03}. 
Let 
$
\A = {\scriptsize \begin{pmatrix} 0&0&0 \\ 0&-1&0 \\ 0&0&0\end{pmatrix}}, \
\B = {\scriptsize \begin{pmatrix} 1&0&0 \\ 0&0&0 \\ 0&0&1\end{pmatrix}}.   
$
Then 
$
\B_\leq = \{(u_1,u_2) : u_1^2 + 1 \leq 0 \} = \emptyset \subseteq   
\A_\geq = \{(u_1,u_2) : -u_2^2 \geq 0\} = \{(u_1,0)\}. 
$ 
But $\coneJ_-(\B) \ni {\scriptsize \begin{pmatrix} 0&0&0 \\ 0&1&0 \\ 0&0&0\end{pmatrix}}
\not\in \coneJ_+(\A)$; hence $\coneJ_-(\B) \not\subseteq \coneJ_+(\A)$. 
\erema
%%%%%%%%%

\medskip

Now we are ready to prove Theorem~\ref{theorem:main2}. 
By Lemmas~\ref{lemma:NScond03}, we see that for every $\A, \ \B \in \SymMat^n_+$ 
\begin{eqnarray*}
 \emptyset \not= \B_\leq \subseteq \A_\geq 
& \Rightarrow &\coneJ_0(\B) \subseteq \coneJ_-(\B) \subseteq \coneJ_+(\A), \\ 
 \emptyset \not= \A_\geq = (-\A)_\leq \subseteq \B_\geq
& \Rightarrow & \coneJ_+(\A) = \coneJ_-(-\A) \subseteq \coneJ_+(\B). 
\end{eqnarray*}
(Note that the assumption $\BC_\geq \not= \emptyset$ of condition (C) with $\coneF = \SymMat^n_+$ implies $\A_\geq \not= \emptyset$ for every $\A \in \BC$.) 
Therefore, condition (C) implies condition (B). 
\qed

\subsection{Proof of Theorem~\ref{theorem:JOYCE2024ext}}

\label{section:extension}

We have shown in [g] that 
condition (D) is equivalent to condition ($\widetilde{\rm D}$).  
By Theorem~\ref{theorem:JOYCE2024}, we obtain 
\begin{eqnarray*}
\begin{array}{l}
\coneJ_+(\widetilde{\BC})\cap\widetilde{\spaceH^1}=
\overline{\mbox{co}}(\coneG^r\cap\coneJ_+(\widetilde{\BC})\cap\widetilde{\spaceH^1}),\\
[3pt]
\eta(\coneJ_+(\widetilde{\BC}),\widetilde{\Q},\widetilde{\H^1}) 
=  \eta(\coneG^r\cap\coneJ_+(\widetilde{\BC}),\widetilde{\Q},\widetilde{\H^1}), 
\end{array}
\end{eqnarray*}
which are equivalent to~\eqref{eq:JOYCE2024a} and~\eqref{eq:JOYCE2024b}, 
where $\widetilde{\spaceH^1} = \{ \V \in \SymMat^r : \inprod{\widetilde{\H^1}}{\V} = 1 \}$.
\qed  

\section{Relationships among conditions (I), (II), (III), (B), (C), and (D)} %  under moderate assumptions} 

\label{section:equivalence}

Throughout this section, we assume 
$\coneF = \coneF_{\min}$ in conditions (B), (C) and (D). 
We have already observed that 
$\mbox{(I)} \Rightarrow \mbox{(B)}, \ \mbox{(III)} \Rightarrow \mbox{(D)} \ \mbox{and } 
\ \mbox{(C)} \Rightarrow \mbox{(B)}$ (Theorem~\ref{theorem:main2}).  
To prove that condition (B) is the weakest of these conditions, 
in the sense it is implied by all the other conditions, 
we show $\mbox{(II)} \Rightarrow \mbox{(B)}$ in Section~\ref{section:IItoB}, 
and $\mbox{(D)} \Rightarrow \mbox{(B)}$ 
in Section~\ref{section:DtoB}.
Under appropriate assumptions, equivalence of 
conditions (I), (II), (B) and (C) are shown in Section~\ref{section:IandIIBC}.

\subsection{Condition (II) $\Rightarrow$ condition (B)}

\label{section:IItoB}

Assume that condition (II) is satisfied. Let $\A, \B \in \BC$ and $\A \not= \B$.  
Then, there exists an $(\alpha,\beta) \not= \0$ such that 
$\alpha\A+\beta\B \in \SymMat^n_+$, which implies 
\begin{eqnarray}
\inprod{\alpha\A+\beta\B}{\X} \geq 0 \ \mbox{for every } \X \in 
 \SymMat^n_+. 
\label{eq:ARGUEineq}
\end{eqnarray}
We will show that $\coneJ_0(\B) \cap \coneF \subseteq \coneJ_+(\A) \cap \coneF$ or 
$\coneJ_+(\A) \cap \coneF \subseteq \coneJ_+(\B) \cap \coneF$ holds. 
One of the following cases (a),(b),$\ldots,$(f) occurs. \vspace{-2mm}
\begin{description}
\item{(a) } 
$\alpha > 0$:  
If $\X \in \coneJ_0(\B) \cap \coneF$, then we see from~\eqref{eq:ARGUEineq} 
that $\X \in \coneJ_+(\A) \cap \coneF$. Hence, 
$\coneJ_0(\B) \cap \coneF \subseteq \coneJ_+(\A) \cap \coneF$ holds.  \vspace{-2mm}
\item{(b) }  $\alpha = 0$ and $\beta > 0$: By \eqref{eq:ARGUEineq}, 
$\coneJ_+(\A)\cap\coneF \subseteq \coneF = \coneJ_+(\B)\cap \coneF $ holds. 
\vspace{-2mm}
\item{(c) } $\alpha = 0$ and $\beta < 0$: By \eqref{eq:ARGUEineq}, 
$\inprod{-\B}{\X} \geq 0$ holds for every  $\X \in \SymMat^n_+$, which implies 
$-\B \in \SymMat^n_+$ and $\coneJ_+(\B)$ forms a face of $\SymMat^n_+$. 
If $\coneJ_+(\B) \cap \coneF = \coneF$, then 
$\coneJ_+(\A) \cap \coneF \subseteq \coneJ_+(\B) \cap \coneF$ holds. Otherwise, 
$\coneJ_+(\BC) \cap \coneF$ 
is contained in the proper face $\coneJ_+(\B) \cap \coneF$ 
of $\coneF = \coneF_{\min}$, which contradicts the definition of $\coneF_{\min}$. 
\vspace{-2mm}
\item{(d) } $\alpha < 0$ and $\beta > 0$: 
Then $\B = (-\alpha/\beta)\A+\Y $ for some $\Y \in \SymMat^n_+$, 
which implies $\coneJ_+(\A) \cap \coneF \subseteq \coneJ_+(\B) \cap \coneF$. 
\vspace{-2mm}
\item{(e) } 
$\alpha < 0$, $\beta =0$. This case can be  treated similarly to case (c) to show that 
$\coneJ_+(\A) \cap \coneF$ forms a face of $\coneF = \coneF_{\min}$. 
\vspace{-2mm}
\item{(f) } 
$\alpha < 0$, $\beta < 0$. In this case, we observe that 
\begin{eqnarray*}
\coneJ_+(\BC) & \subseteq  & \left\{ \X \in \coneF :  \inprod{\A }{ \X} \geq 0, \ \inprod{\B}{ \X} \geq 0 \right\} \\
& = & \left\{ \X \in\coneF  :  \inprod{\A}{ \X} \geq 0, \ \inprod{\B}{ \X} \geq 0, \  -\inprod{\alpha\A+ \beta\B}{ \X} \geq 0 \right\}\\
& & \mbox{(since $\alpha < 0$ and $\beta <0$)} \\
% & & mbox{(since $-(\alpha\A+ \beta\B) \in \SymMat^n_+$)} \\
& \subseteq & \left\{ \X \in \coneF  : \inprod{-\alpha\A - \beta \B}{ \X} = 0 \right\} \
 \mbox{(since $\alpha\A + \beta \B \in \SymMat^n_+$)}. 
\end{eqnarray*}
Hence, $\coneJ_+(\BC)$ is included in a face 
$\coneF' \equiv \{\X \in \coneF : \inprod{-\alpha\A -\beta \B}{ \X} = 0\}$ of $\coneF$.  
If $\coneF' = \coneF$ then 
\begin{eqnarray*}
\coneJ_0(\B) \cap \coneF = \coneJ_0(\B) \cap \coneF' 
= \coneJ_0(\A) \cap \coneF' = \coneJ_0(\A) \cap \coneF \subseteq \coneJ_+(\A) \cap \coneF. 
\end{eqnarray*}
Otherwise, $\coneF'$ is a proper face of $\coneF = \coneF_{\min}$ that contains 
$\coneJ_+(\BC)$, a contradiction to the definition of $\coneF_{\min}$. 
\vspace{-4mm} 
%\vspace{-2mm}
\end{description}
\qed 

\subsection{Condition (D) $\Rightarrow$ condition (B)}

\label{section:DtoB}

We have shown the equivalence of conditions (B) and ($\widetilde{\rm B}$) 
in Section~\ref{section:facialReductionOfJ+} [c] and the equivalence of 
(D) and ($\widetilde{\rm D}$) in Section~\ref{section:facialFBmain2} [g]. 
Therefore, it suffices to show 
($\widetilde{\rm D}$)  $\Rightarrow$ ($\widetilde{\rm B}$)
 for  (D) $\Rightarrow$  (B). 
 % }. 
Since we take $\coneF = \coneF_{\min}$, 
$\coneJ_+(\widetilde{\BC}) \cap \SymMat^r_{++} \not= \emptyset$ 
(Slater's constraint qualification for SDP~\eqref{eq:SDP10}). 
For simplicity of notation, we omit $\widetilde{ \ }$ and assume that 
$\coneJ_+(\BC)$ itself satisfies condition \vspace{-2mm}
\begin{description}
\item{(A-1) } $\coneJ_+(\BC) \cap \SymMat^n_{++} \not= \emptyset$.\vspace{-2mm}
\end{description}
throughout this and the next sections.  

\medskip

For the proof of  (D) $\Rightarrow$ (B), 
we fix arbitrary $\A, \ \B \in \BC$. 
The case $\B\in\SymMat^n_+$ is immediate. Indeed, if $\B\in\SymMat^n_+$, then
$\inprod{\B}{\X} \geq 0$ for every $\X \in \SymMat^n_+$,  
and hence $\coneJ_+(\B) = \SymMat^n_+$. 
Therefore $\coneJ_+(\A)\subseteq\coneJ_+(\B)$, and condition (B) holds.
Thus we assume $\B\not\in\SymMat^n_+$. 
Let us write % $\A, \ \B  \in \SymMat^n$ as 
$ 
\A = {\scriptsize \begin{pmatrix}\C& \cc \\ \cc^T & \gamma \end{pmatrix}}, \
\B = {\scriptsize \begin{pmatrix}\D& \d \\ \d^T & \delta \end{pmatrix}},
$ 
where $\C, \D \in \SymMat^{n-1}$, $\cc, \d \in \Real^{n-1}$, and $\gamma, \delta \in \Real$.
% where $\C, \ \D \in\SymMat^{n-1}, \ \cc, \ \d \in \Real^{n-1}, \  \gamma, \ \delta \in \Real$. 
We show that $\B_= \subseteq \A_\geq$ $\Rightarrow$ $\coneJ_0(\B) \subseteq \coneJ_+(\A)$ 
under the assumptions that
\begin{eqnarray}
& &  \coneJ_+(\B) \cap \SymMat^n_{++} \not= \emptyset, \ \B \not\in \SymMat^n_+, \ \mbox{and } \D \not= \O. 
\label{eq:condition0}
\end{eqnarray}
Define the quadratic functions $h$ and $f$ on $\Real^{n-1}$ as 
\begin{eqnarray*}
h(\u) = \u^T\D\u + 2\d^T\u + \delta, \ f(\u) = \u^T\C\u + 2\cc^T\u + \gamma \ \mbox{for every } \u \in \Real^{n-1}. 
\end{eqnarray*}

%%%%%%%%%
\lemm \label{lemma:NScond02} \mbox{ \ } \vspace{-2mm} 
\begin{description}
\item{(i) } The quadratic function $h : \Real^{n-1} \rightarrow \Real$ takes both negative and positive values; 
$h(\u^1) < 0 < h(\u^2)$ for some $\u^1, \ \u^2 \in \Real^{n-1}$.\vspace{-2mm} 
\item{(ii) } The following two conditions are equivalent.  \vspace{-2mm} 
\begin{description}
\item{$\mbox{E}_1$: } $h(\u) = 0$ $\Rightarrow$ $f(\u) \geq 0$. \vspace{-1mm}  
\item{$\mbox{E}_2$: } There exists $\tau \in \Real$ such that $f(\u) + \tau h(\u)$ for every $\u \in \Real^{n-1}$.\vspace{-2mm}  
\end{description}
\end{description}
\elemm
%%%%%%%%%
\proof{
(i) The existence of $\u^1$ follows from $\B \not\in \SymMat^{n}_+$, and that of $\u^2$ follows from 
$\coneJ_+(\B) \cap \SymMat^n_{++} \not= \emptyset$. (ii) Under the assumptions (i) and $\D \not= \O$, 
the equivalence of $\mbox{E}_1$ and $\mbox{E}_2$ follows from \cite{XIA2016}, which established 
the so-called S-Lemma with equality, and provided the equivalence of $\mbox{E}_1$ and $\mbox{E}_2$ under 
appropriate assumptions. In particular, Theorem 3 of \cite{XIA2016} stated a necessary and sufficient 
condition for the equivalence 
% of $\mbox{E}_1$ and $\mbox{E}_2$ 
when (i) holds, where $\D \not= \O$ was shown to be sufficient for the 
equivalence. 
\qed

\medskip

Obviously, $\B_= \subseteq  \A_\geq$ is equivalent to $\mbox{E}_1$. 
Hence, by Lemma~\ref{lemma:NScond02}, there exists a $\tau \in \Real$ such that $f(\u) + \tau h(\u)$ for every $\u \in \Real^{n-1}$.  
By \cite[Lemma 3.3]{LAURENT2008}, this is equivalent to $\A + \tau\B \in \SymMat^n_+$, that is, 
$\inprod{\A + \tau\B}{\X} \geq 0$  for every $\X \in \SymMat^n_+$. 
Consequently, if $\inprod{\B}{\X}= 0$ with $\X \in \SymMat^n_+$,  
then $\inprod{\A}{\X} \geq 0$. This shows that $\coneJ_0(\B) \subseteq \coneJ_+(\A)$.
\qed
}

\subsection{Equivalence of conditions (I), (II), (B), and (C)}

\label{section:IandIIBC}

Throughout this section, we allow
 $\BC \subseteq \SymMat^n$ to be infinite  
in conditions (I) and (II) as in condition (B)  
even though conditions (I) and (II) are originally stated for finite $\BC \subseteq \SymMat^n$.
In addition to $\coneF = \SymMat^n_+$ and condition (A-1) assumed in the previous section, 
we assume conditions 
\vspace{-1mm}
\begin{description}
\item{(A-2) } $\coneJ_+(\B) \not\subseteq \coneJ_+(\A)$ for every distinct $\A, \B \in \BC$.  
 \vspace{-2mm}
 \item{(A-3) } $\B_\leq \not= \emptyset$ for every $\B \in \BC$. 
 \vspace{-1mm}
\end{description}
Under these conditions, we show the equivalence of conditions (I), (II), (B), and (C). 
In the case where $\BC$ is finite, 
if 
$\coneJ_+(\A) \subseteq \coneJ_+(\B)$ for distinct $\A, \ \B \in \BC$,  
we can remove $\B \in \BC$ from $\BC$ one by one, recursively,  so that 
$\coneJ_+(\BC) =  \coneJ_+(\BC\backslash\{\B\})$. Therefore, condition (A-2) is 
attained. 
When $\BC$ is infinite, 
we need Zorn's lemma (see, for example, \cite{HRBACEK1999,JECK2006}) 
to remove such $\B$'s consistently from $\BC$ so that the resulting $\BC$ 
satisfies condition (A-2). The details are omitted here. Condition (A-3) is also reasonable 
since if $\B_\leq = \emptyset$ then $\B_\geq = \Real^{n-1}$ and $\BC_\geq = 
\left\{\BC\backslash\{\B\}\right\}_\geq$. 

\medskip

The equivalence of (I) and (B) is obvious under (A-2). 
We have already seen that (II) $\Rightarrow$ (B) in Section~\ref{section:IItoB}  and that 
(C) $\Rightarrow$ (B) in Theorem~\ref{theorem:main2}. To prove 
(B) $\Rightarrow$ (II) and (B) $\Rightarrow$ (C), it suffices to prove the following lemma, 
which also shows that (B) implies~\eqref{eq:condJOYCE2024ext} included in (D). 

\medskip

%%%%%%%%%
\lemm \label{lemma:equivalence0}
Assume that $\BC \subseteq \SymMat^n$ satisfies conditions (A-1) and (A-2). 
Let $\A, \B \in \BC$ and $\A \not= \B$. Then 
 \begin{eqnarray}
& & \alpha\A + \beta\B \in \SymMat^n_+ \ \mbox{for some } (\alpha,\beta) \not= \0.  \label{eq:equivalenceC0}\\
& & \hspace{13mm} \Uparrow \nonumber \\
& &   \coneJ_0(\B) \subseteq \coneJ_+(\A) \ \Leftrightarrow \  \coneJ_-(\B) \subseteq \coneJ_+(\A) 
\label{eq:equivalenceA0}\\
& & \hspace{13mm} \Downarrow \hspace{35mm} \Downarrow 
\nonumber \\ 
& & \hspace{5mm} \B_= \subseteq \A_\geq  \hspace{20mm} 
\B_\leq \subseteq \A_\geq.  \nonumber 
 \end{eqnarray}
\elemm
%%%%%%%%%
\proof{
\begin{description}
\item{(i) } In~\eqref{eq:equivalenceA0}, $\Leftarrow$  is straightforward since 
$\coneJ_0(\B) \subseteq \coneJ_-(\B)$. To prove $\Rightarrow$, 
assume on the contrary 
that $\coneJ_0(\B) \subseteq \coneJ_+(\A)$ but 
$\coneJ_-(\B) \not\subseteq \coneJ_+(\A)$  or equivalently that 
$\inprod{\B}{\overline{\X}} < 0$ \mbox{and } 
$\inprod{\A}{\overline{\X}}<0$ for some $\overline{\X} \in \SymMat^n_+$. 
By condition (A-2), 
$\coneJ_+(\B) \not\subseteq \coneJ_+(\A)$, which together with 
$\coneJ_0(\B) \subseteq \coneJ_+(\A)$ implies 
%} 
% 
$\inprod{\B}{\widetilde{\X}} > 0$ \mbox{and } 
$\inprod{\A}{\widetilde{\X}}<0$ for some $\widetilde{\X} \in \SymMat^n_+$.
Hence there exists $\lambda \in (0,1)$ such that 
  \begin{eqnarray*}
  \inprod{\B}{\lambda\overline{\X}+ (1-\lambda)\widetilde{\X}} = 0, \  
  \inprod{\A}{\lambda\overline{\X}+ (1-\lambda)\widetilde{\X}} <  0, \ \lambda\overline{\X}+ (1-\lambda)\widetilde{\X} \in \SymMat^n_+.  
  \end{eqnarray*} 
  This contradicts  $\coneJ_0(\B) \subseteq \coneJ_+(\A)$. \vspace{-2mm} 
\item{(ii) } $\coneJ_0(\B) \subseteq \coneJ_+(\A) \  \Rightarrow \  \eqref{eq:equivalenceC0}$: 
Consider the primal-dual pair of  SDPs 
\begin{eqnarray}
\zeta_{\rm p}  & = & \inf\{ \inprod{\A}{\X} : \X \in \SymMat^n_+, \ \inprod{\B}{\X} = 0 \} = \inf\{ \inprod{\A}{\X} : \X \in \coneJ_0(\B)\}, \label{eq:primalSDP} \\
\zeta_{\rm d}  & = & \sup\{ 0 : \A + \tau \B \in \SymMat^n_+, \ \tau \in \Real \}. \nonumber
\end{eqnarray}
Obviously, $\zeta_{\rm p} = 0$ if and only if  $\coneJ_0(\B) \subseteq \coneJ_+(\A)$. 
By (A-1),  $-\B \not\in \SymMat^n_+$, which implies $\inprod{\B}{\X^1} > 0$ for some 
$\X^1 \in \SymMat^n_{++}$. By (A-2), $\B \not\in \SymMat^n_+$, which implies 
$\inprod{\B}{\X^2} < 0$ for some $\X^2 \in \SymMat^n_{++}$. 
Hence a convex combination 
$\overline{\X}\in\SymMat^n_{++}$ of $\X^1$ and $\X^2$ satisfies $\inprod{\B}{\overline{\X}} = 0$, {\it i.e.}, $\overline{\X}$ is an interior feasible solution of 
primal SDP~\eqref{eq:primalSDP}. By the duality theorem,  $\zeta_p = \zeta_d = 0$ 
if and only if the dual SDP is feasible, {\it i.e.}, $\A + \tau \B \in \SymMat^n_+$ for some 
$\tau \in \Real$. Therefore, we have shown that 
$\coneJ_0(\B) \subseteq \coneJ_+(\A) \  \Rightarrow \  \eqref{eq:equivalenceC0}$. %  has been shown.
\vspace{-2mm}
\item{(iii) } $\coneJ_-(\B) \subseteq \coneJ_+(\A)$ $\Rightarrow$ $\B_\leq \subseteq \A_\geq$: 
We observe that 
\begin{eqnarray*}
& \coneJ_-(\B) \subseteq \coneJ_+(\A)& \\
& \hspace{38mm} \Updownarrow \ \mbox{(by Lemma~\ref{lemma:NScond00} (ii))} & \\
&  \{\x \in \Real^n : \inprod{\B}{\x\x^T} \leq 0 \} 
\subseteq \{\x \in \Real^n : \inprod{\A}{\x\x^T} \geq 0 \} & \\ 
& \Downarrow & \\
&  \{\x \in \Real^n : \inprod{\B}{\x\x^T} \leq 0, \ x_n = 1 \} 
\subseteq \{\x \in \Real^n : \inprod{\A}{\x\x^T} \geq 0, \ x_n = 1\} & \\ 
& \Updownarrow & \\
&  \B_\leq \subseteq \A_\geq. &
\end{eqnarray*}
\vspace{-10mm}
\item{(iv) } $\coneJ_0(\B) \subseteq \coneJ_+(\A)$ $\Rightarrow$ $\B_= \subseteq \A_\geq$: 
This assertion can be proved similarly as in (iii).
%  by applying Lemma~\ref{lemma:NScond02} instead of Lemma~\ref{lemma:NScond00} (ii).
\end{description}
\qed
}

\rema \label{remark:dualOfB}
From (ii) of the proof above, we see that condition (II) can be regarded as the dual of 
condition (B) under assumptions (A-1) and (A-2). 
\erema

\section{Examples}

\label{section:examples}

In \cite[Section 4.1]{ARIMA2023},
several examples satisfying condition (B) 
with finite $\BC$ and $\coneF = \SymMat^n_+$ were provided.    
We  present  three examples that are not covered by those examples in this section.
%

%%%%%%%%%
\examp \label{example:infiniteBC}
This example provides an infinite $\BC \subseteq \SymMat^n$ satisfying condition (C). 
Let 
\begin{eqnarray*}
& & \B(\t) = \begin{pmatrix} \I & - \t \\ -\t^T & \t^T\t - r^2 \end{pmatrix} \in \SymMat^n \ 
(\t \in T),  \ \BC = \left\{\B(\t)  : \t \in T \right\}, 
\end{eqnarray*}
where $0 < r \leq 1/2$,  
$T \subseteq \Integer^n$ (the set of integer vectors in $\Real^n$) and $\I$ denotes the 
$(n-1) \times (n-1)$ identity matrix.  Then, 
\begin{eqnarray*}
& & q(\u,\B(\t)) = \begin{pmatrix} \u\\ 1\end{pmatrix}^T \B(\t)\begin{pmatrix} \u\\ 1\end{pmatrix} = 
\parallel \u - \t \parallel^2- r^2, \\
& & \B(\t)_\geq \ \mbox{or } \B(\t)_\leq = 
\left\{ \u \in \Real^{n-1} :  \parallel \u - \t  \parallel ^2- r^2  \geq 0 \ \mbox{or  $\leq 0$, respectively}
\right\}
\end{eqnarray*} 
for every  $\t \in T$ and $\u \in \Real^{n-1}$.  See Figure 3. 
It is easily seen that condition (C) with $\coneF = \SymMat^n_+$ is satisfied.
Therefore, by Theorem~\ref{theorem:main2}, we obtain $\coneJ_+(\BC) \in \wFC(\coneG^n)$. 
When $T$ is finite, this example is a special case of quadratic programs with hollows 
\cite{YANG2018}.

\medskip

As a generalization, %we can easily
it is straightforward to construct an ellipsoid-based constraint by replacing each $\B(\t)$ with $\L^T\B(\t)\L$ 
$(\t \in T)$, 
where $\L$ denotes an $n \times n$ nonsingular matrix of the form $\L = \begin{pmatrix} \M & \0 \\ \0^T & 1 \end{pmatrix}$. 
We also note that the equivalence relation~\eqref{eq:wFC} 
between COP($\coneJ\cap\coneG^n,\Q,\H$) 
and its SDP relaxation COP($\coneJ,\Q,\H$) with $\coneJ =\coneJ_+(\BC)$ (or  
$\coneJ =\coneJ_+(\{\L^T\B\L : \B \in \BC\})$)
holds for any choice of $\Q \in \SymMat^n$ and $\H \in \SymMat^n$ 
by Theorem~\ref{theorem:propFC}. 
For example, we can take  $\H = \delta \I$  for some $\delta  > 0$ where $\I$ is the $n \times n$ identity matrix. 
In this case, COP($\coneJ\cap\coneG^n,\Q,\H$) turns out to be 
\begin{eqnarray*}
\eta(\coneJ\cap\coneG^n,\Q,\H) & = & 
\inf\left\{\begin{pmatrix} \u \\ \z \end{pmatrix}^T\Q\begin{pmatrix} \u \\ \z \end{pmatrix} :
\begin{array}{l}
\parallel \u - \t z \parallel ^2- r^2 z^2 \geq 0 \ (\t \in T),  \\[3pt] 
\mbox{(or $\parallel \M\u - \t z \parallel ^2- r^2 z^2 \geq 0 \ (\t \in T)$)}, \\[3pt] 
u_1^2+ \cdots + u_{n-1}^2 + z^2 = 1/\delta
\end{array}
\right\}. 
\end{eqnarray*}
\eexamp
%%%%%%%%%

\begin{figure}[t!]  \vspace{-0.3cm} 
\begin{center}
\includegraphics[height=50mm]{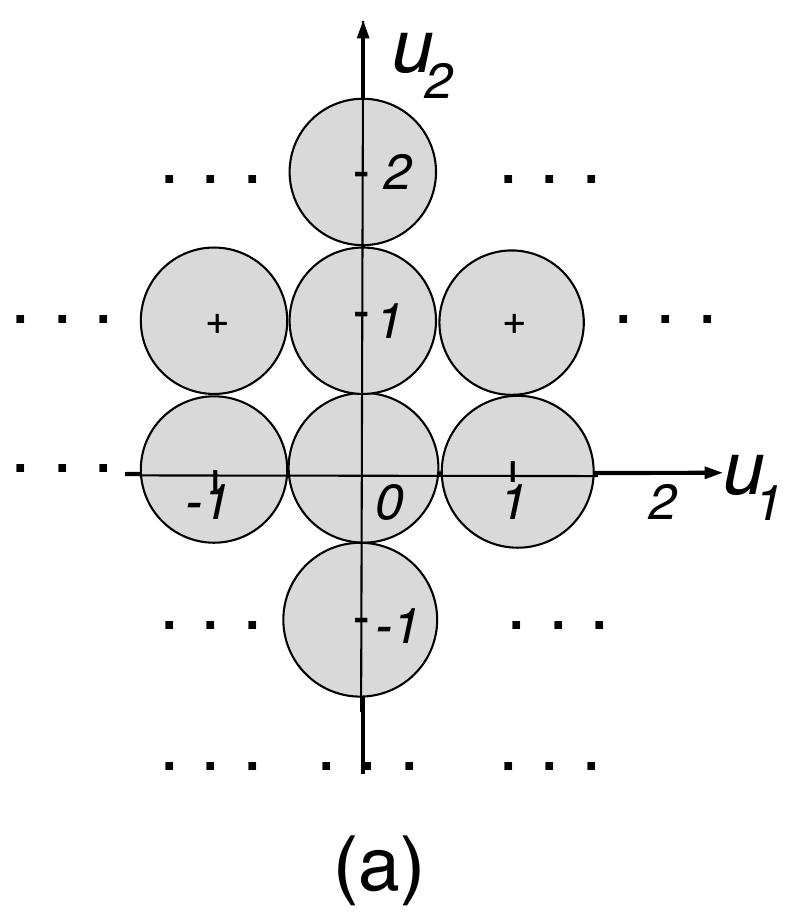}
\hspace{20mm} 
\includegraphics[height=50mm]{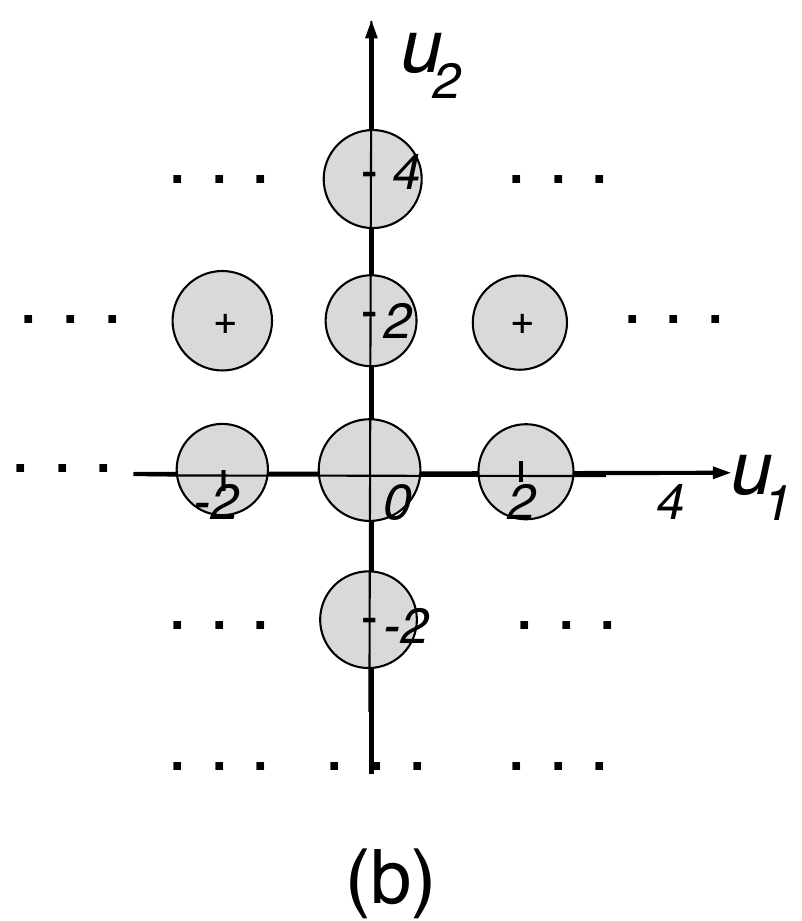}
\end{center}

\vspace{-5mm}

\caption{Illustration of $\B(\t)_\leq$. Each  gray disk region
corresponds to $\B(\t)_\leq$ for some $\t \in T$. (a) $T = \Integer^{n-1}$ and $r = 1/2$. 
(b) $T = 2  \Integer^{n-1}$ and $r = 1/3$. 
 \label{fig:Fig3}
}
\end{figure} 

\medskip

%%%%%%%%%
\examp \label{example:extensionH}  
This example presents another infinite $\BC$ satisfying condition (C).
For  a finite or infinite sequence $\{(\gamma_p,\mu_p,\lambda_p,\sigma_p) \in \Real^4:p=1,2,\ldots\}$ satisfying 
\begin{eqnarray}
\gamma_p > 0,\ \mu_p \geq 0 ,\lambda_p > 0 \ \mbox{for every } p, \ \mbox{and } 
|\sigma_q - \sigma_p| \geq \sqrt{\lambda_p} + \sqrt{\lambda_q} \ \mbox{if } p < q,   \label{eq:parameter}
\end{eqnarray}
we consider a sequence of quadratic functions 
\begin{eqnarray}
q_p(\u,\v,\w) & = & 
\parallel \v - \sigma_p \u \parallel^2 - \lambda_p\parallel \u\parallel^2 + \mu_p\parallel\w\parallel^2 + \gamma_p \nonumber \\
& \ & \mbox{for every } \u \in \Real^\ell, \ \v \in \Real^\ell   \ \mbox{and } 
\w \in \Real^m  \label{eq:quadraticFunction}
\end{eqnarray}
$(p=1,2,\ldots)$. 
We can take  a matrix $\B^p \in \SymMat^{2\ell+m+1}$  representing the quadratic function 
$q^\sigma(\u,\v,\w)$ such that 
\begin{eqnarray*}
q^\sigma(\u,\v,\w) 
& = &
{\scriptsize\begin{pmatrix}\u\\ \v \\ \w \\ 1\end{pmatrix}}^T\B^p{\scriptsize\begin{pmatrix}\u\\ \v \\ \w \\ 1\end{pmatrix}}
\ \mbox{for every } \u \in \Real^\ell, \ \v \in \Real^\ell   \ \mbox{and } 
\w \in \Real^m     
\end{eqnarray*}
$(p=1,2,\ldots)$. 

\medskip

We show that $\BC = \{\B^p : p=1,2,\ldots\}$ satisfies condition (C) with 
$\coneF = \SymMat^{2\ell+m+1}_+$. Let $p < q$. 
Assume on the contrary that $(\u,\v,\w) \in \B^p_\leq \cap \B^q_\leq$ 
for some $(\u,\v,\w) \in \Real^{2\ell+m}$. Then 
\begin{eqnarray*}
\parallel \v - \sigma_p \u \parallel^2 + \mu_p\parallel \w\parallel^2 + \gamma_p \leq \lambda_p\parallel \u\parallel^2, \\  
\parallel \v - \sigma_q \u \parallel^2 + \mu_q\parallel \w\parallel^2 + \gamma_q \leq \lambda_q\parallel \u\parallel^2.  
\end{eqnarray*}
Hence,
\begin{eqnarray*}
& & 0 < \parallel \u \parallel, \ \parallel \v - \sigma_p \u \parallel <  \sqrt{\lambda_p}\parallel \u\parallel, \ 
\parallel \v - \sigma_q \u \parallel <  \sqrt{\lambda_q}\parallel \u\parallel.  
\end{eqnarray*}
Therefore,
\begin{eqnarray*}
|\sigma_q-\sigma_p| \parallel \u \parallel & = & \parallel (\sigma_q \u-\v) + (\v - \sigma_p \u)\parallel \\
& \leq &
\parallel\sigma_q \u-\v\parallel + \parallel\v - \sigma_p \u\parallel \\
& < & (\sqrt{\lambda_p}+\sqrt{\lambda_q})\parallel \u\parallel, 
\end{eqnarray*}
which implies $\sigma_q-\sigma_p < \sqrt{\lambda_p}+\sqrt{\lambda_q}$, a contradiction to 
the last inequality of \eqref{eq:parameter}. Thus we have shown 
$\B^p_\leq \cap \B^q_\leq  = \emptyset$, which implies  $\B^p_\leq \subseteq \B^q_\geq$ 
and $\B^q_\leq \subseteq \B^p_\geq$. 

\medskip

If we take $\gamma_p = \lambda_p = \mu_p = 1$ and $\sigma_p = 2p$ $(p=1,2,\ldots)$, then 
$\{(\gamma_p,\mu_p,\lambda_p,\sigma_p):p=1,2,\ldots\}$ satisfies \eqref{eq:parameter}. 
As another case, we consider 
\begin{eqnarray}
\gamma_p > 0, \ \lambda_p = (a_{p}-a_{p-1})^2/4, \ \mu_p = 0, \ \sigma_p = -(a_{p-1}+a_p)/2 \ (p=1,2,\ldots), \label{eq:specialCase}
\end{eqnarray}
where $\{a_p:p=0,1,\ldots \}$ denotes an infinite sequence of positive numbers such that 
$a_{p-1} < a_p$ In this case,  if $p < q$ then 
\begin{eqnarray*}
|\sigma_q - \sigma_p| - (\sqrt{\lambda_q} + \sqrt{\lambda_p}) & = &
\frac{(a_q+a_{q-1})-(a_p+a_{p-1})}{2} - \frac{(a_q-a_{q-1}) + (a_p-a_{p-1})}{2}\\
 & = & a_{q-1} -a_p \geq 0. 
\end{eqnarray*}
Therefore, the sequence $\{(\gamma_p,\lambda_p,\mu_p,\sigma_p) : p=1,2,\ldots\}$ defined by 
\eqref{eq:specialCase} satisfies \eqref{eq:parameter}. Now suppose that $\ell=1$ and $m=0$. 
Then, the quadratic functions \eqref{eq:quadraticFunction} turns out to be 
\begin{eqnarray*}
q_p(u,v) & = & (v-\sigma_pu)^2  - \lambda_pu^2 + \gamma_p \\
& = & (v-(a_{p-1}+a_p)u/2)^2  - (a_{p}-a_{p-1})^2u^2/4  + \gamma_p \\ 
% 
%& = & \TBrown{v^2- (a_{p-1}+a_p)uv + (a_{p-1}+a_p)^2u^2/4- (a_{p}-a_{p-1})^2u^2/4  + \gamma_p} \\
%& = &  \TBrown{v^2- (a_{p-1}+a_p)uv + a_{p-1}a_pu^2 + \gamma_p} \\
%
& = & a_{p-1}a_pu^2 - (a_{p-1}+a_p)uv + v^2 + \gamma_p \\
%& = &  (v-a_{p-1}u)(v - a_pu) + \gamma_p^2. % \\
& = & {\scriptsize \begin{pmatrix} u \\ v \\ 1 \end{pmatrix}}^T \B^p {\scriptsize \begin{pmatrix} u \\ v \\ 1 \end{pmatrix}}, 
\end{eqnarray*}
where 
\begin{eqnarray*}
\B^p & = & \begin{pmatrix} a_{p-1}a_p & -(a_{p-1}+a_p)/2 & 0 \\ -(a_{p-1}+a_p)/2 & 1 & 0 \\ 
0 & 0 & \gamma_p \end{pmatrix}
\end{eqnarray*}
$\BC = \{\B^p : p=1,2,\ldots\}$ satisfies condition (C) with $\coneF = \SymMat^3_+$. 
See Figure 4. 

Now assume that $a_p \rightarrow \bar{a}$ as $p \rightarrow \infty$. 
Then $\B^p \rightarrow 
 \overline{\B}={\scriptsize \begin{pmatrix}  \bar{a}^2 & -\bar{a}  & \0 \\  -\bar{a} & 1 & 0 \\  0 & 0 & 1 \end{pmatrix}}\not\in \BC$ as $p \rightarrow \infty$. 
 Therefore, $\BC$ is not closed. 
 We also see $\overline{\B} \in \SymMat^3_+$. 
Hence, $\coneJ_+(\B^p) \subseteq \SymMat^3_+ = \coneJ_+(\overline{\B})$ 
$(p=1,2,\ldots,)$.  
Therefore,  $\mbox{cl}\BC$ (the closure of $\BC$) does not satisfy condition (A-2) 
although $\coneJ_+(\mbox{cl}\BC) \in \wFC(\coneG^3)$ remains true. 
\eexamp
%%%%%%%%%

\begin{figure}[t!]  \vspace{-40mm} 
\begin{center}
\includegraphics[height=150mm]{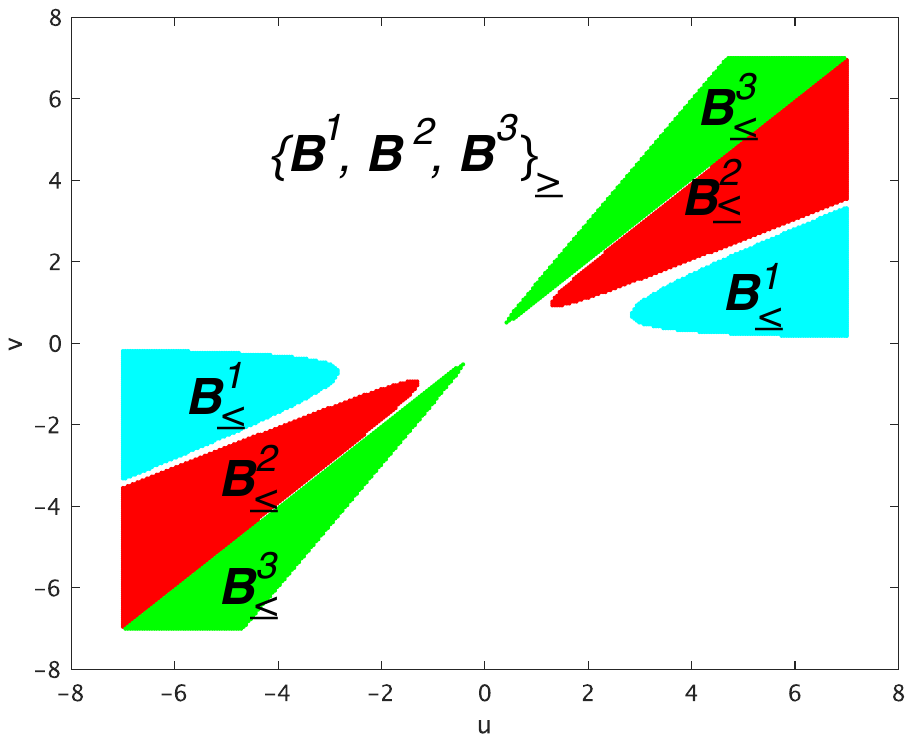}

\vspace{-40mm} 
\caption{
Example~\ref{example:extensionH}.  
We take $a_0 = 0,  \ a_1 = 0.5, \ a_2 = 1, \ a_3 = 1.5, \ \gamma_1=0.5, \ \gamma_2 = 0.1, \ \gamma_2 = 0.01$. 
\label{fig:Fig4}
} 
\end{center}
\end{figure}

\begin{figure}[t!]  \vspace{-30mm} 
\begin{center}
\includegraphics[height=123mm]{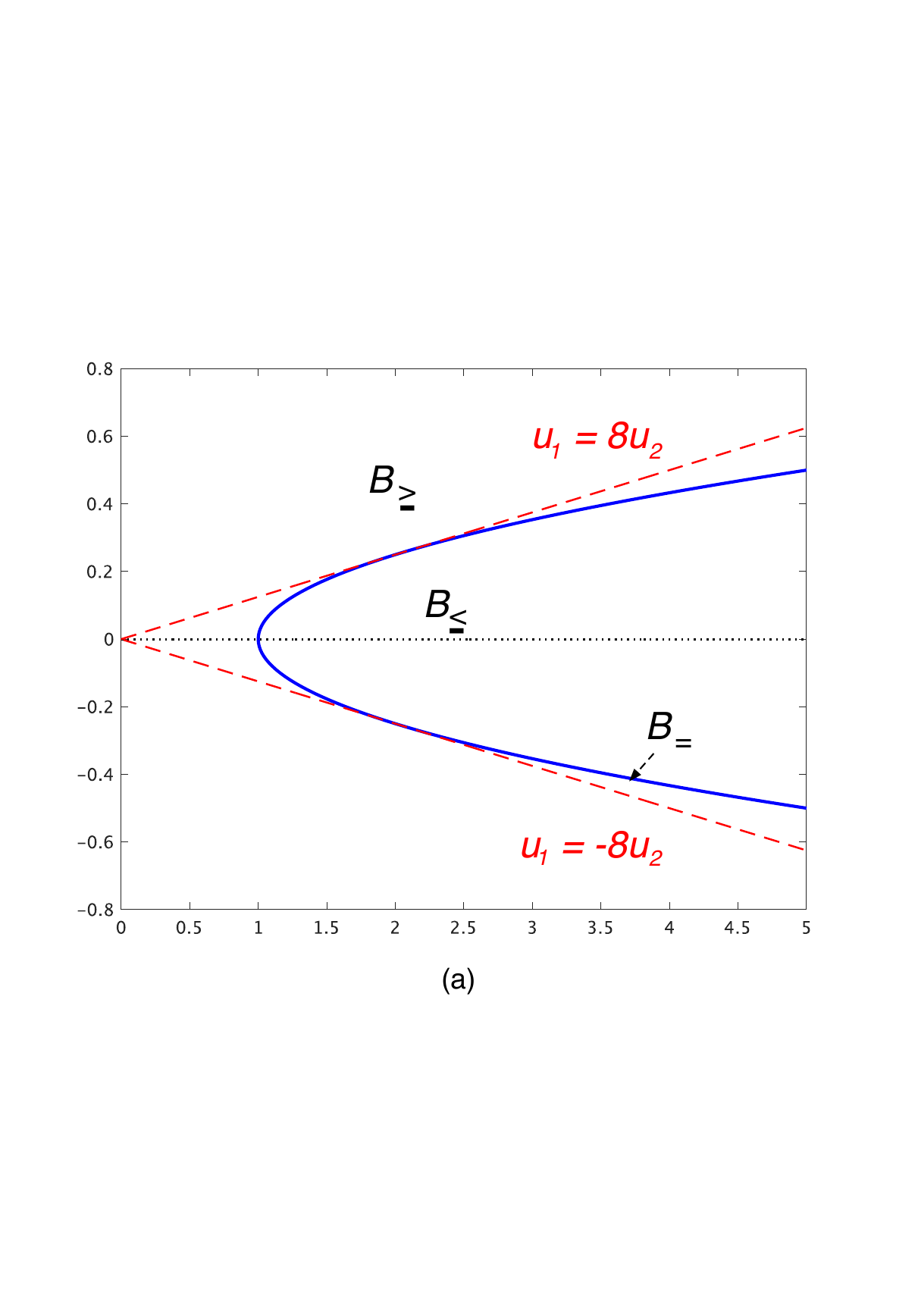}
\hspace{-15mm} 
\includegraphics[height=125mm]{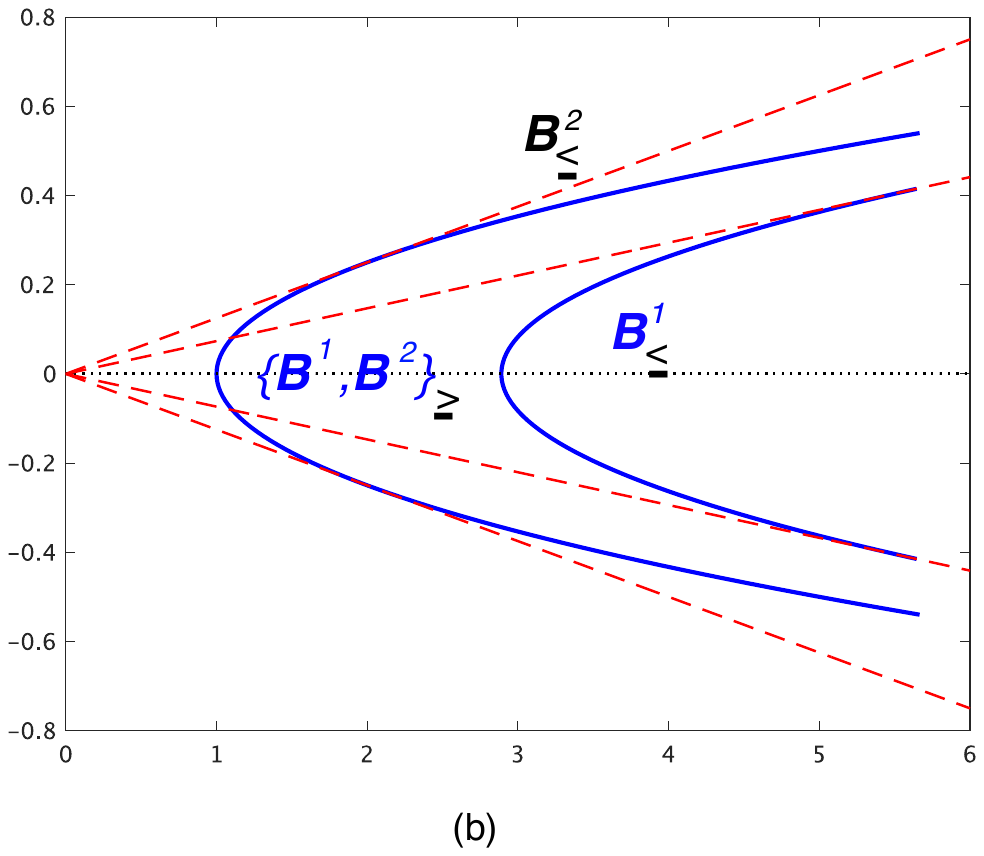}

\vspace{-60mm}

\includegraphics[height=123mm]{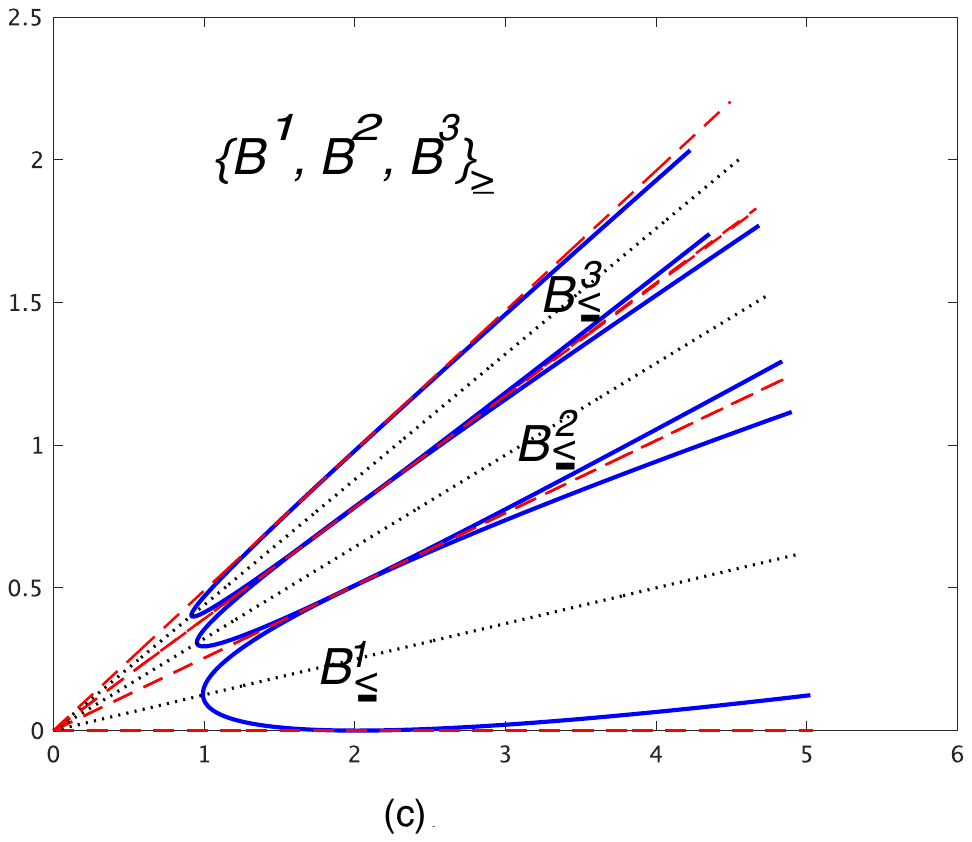}

\end{center}

\vspace{-35mm}

\caption{
Parabola-based constraints with $n= 3$. 
 (a) Parabola $\B_\leq$ defined by~\eqref{eq:parabola} where 
$\lambda_2=16$ and $\lambda_3=1$. 
(b) 
%$\ff_+(1,\{\B_1,\B_2\}) = \ff_+(1,\B_2)\backslash\ff_{--}(1,\B_1)$, where 
% $\B^1_\leq = \{\u \in \Real^2: -u_1 +   16u_2^2 + 3 \leq 0 \}$, 
$\B^1_\geq = \{\u \in \Real^2: -u_1 +   16u_2^2 + 3 \geq 0 \}$, 
%$\B^2_\leq = \{\u \in \Real^2:  -(-u_1 + 16u_2^2 + 1) \leq 0 \}$
%and 
$\B^2_\geq = \{\u \in \Real^2:  -(-u_1 + 16u_2^2 + 1) \geq 0 \}$ 
and $\{\B^1,\B^2\}_\geq = \B^1_\geq \cap \B^2_\geq$. 
%}
\label{fig:Fig5}
}
\end{figure} 

%%%%%%%%%
\examp \label{example:parabola} 
We consider a parabola-based constraint.   Let $n \geq 3$ and 
\begin{eqnarray*}
B_{ij} = \left\{
\begin{array}{ll}
\lambda_i > 0 & \mbox{if } 2 \leq i=j \leq n, \\
-0.5   & \mbox{if } (i,j) = (1,n) \ \mbox{or } (i,j) = (n,1), \\
0    &\mbox{otherwise}. 
\end{array}
\right.
\end{eqnarray*}
Then, 
\begin{eqnarray}
& &  q(\u,\B) =  -u_1 + \sum_{i=2}^{n-1} \lambda_i u_i^2 + \lambda_n  
\ \mbox{for every $\u \in \Real^{n-1}$}, \nonumber \\
& & \B_\geq ,  \B_=  \
 \mbox{or $\B_\leq$}
 = 
 \left\{ \u \in \Real^{n- 1} : 
 \begin{array}{l}
 -u_1 + \sum_{i=2}^{n-1} \lambda_i u_i^2 + \lambda_n \geq 0,  =0  \
\mbox{or  } \leq 0, \\
\mbox{respectively}
\end{array}
\right\}, \label{eq:parabola} \\
& & \B_\leq \subseteq 
\{ \u \in \Real^n :  \lambda_iu_i^2  + \lambda_n \leq u_1  \ (i=2,\ldots,n-1)  \} \nonumber \\
 & &  \hspace{8mm} \subseteq \coneK_-(\B) \equiv 
 \{ \u \in \Real^{n-1} :   0 \leq u_1, \ -u_1 \leq 2\sqrt{\lambda_i\lambda_n}u_i \leq u_1 \ (i=2,\ldots,n-1)\}.  \nonumber
\end{eqnarray}
See Figure 6 (a). We note that 
$\coneK_-(\B)$ forms a polyhedral cone in $\Real^{n-1}$,  which converges to the half line 
$\{ \u \in \Real^{n-1} :u_1 \geq 0, \ u_i = 0  \ (i=2,\ldots,n-1)\}$ as all $\lambda_i$ $(i=2,\ldots,n-1)$ tend to $\infty$. 
We see that each $\B_= \cap \{\u \in \Real^{n-1} : u_j = 0 \ (2 \leq j \not= i \leq n-1)\}$ 
forms a $2$-dimensional 
parabola $(i=2,\ldots,n-1)$. 
By applying a linear transformation $\B \rightarrow \L^T\B\L \in \SymMat^n$ 
with a nonsingular $\L$ to $ \B_=$ with 
different $\lambda_i  > 0 \ (i=2,\ldots,n)$, we can create various 
parabolas. Furthermore, we can arrange some of those parabolas such 
that the assoicated $\BC$  satisfies condition (C). % and (C)'. 
See Figure~6~(b) and~(c).  
\eexamp
%%%%%%%%%

%
%\input sect7ConcludingDiscussions.tex
%!TEX root = ./main.tex

\section{Concluding remarks}

\label{section:conclusion}

Finally, we comment on the role of the representation of the closed convex cone $\coneJ \subseteq \SymMat^n_+$ 
in the geometric QCQP, COP($\coneG^n\cap\coneJ$) and its SDP relaxation COP($\coneJ$). % J+(B).
Our sufficient conditions for the equivalence of COP($\coneG^n\cap\coneJ$) and COP($\coneJ$) 
are formulated in terms of a set $\BC$ % representing the cone 
with $\coneJ = \coneJ_+(\BC)$, %J ⊆ Sⁿ₊, 
and are therefore representation-dependent:
different choices of $\BC$ may affect whether conditions (I), (II), (III), (B), (C), or (D) hold,
even though the cone $\coneJ$ itself is unchanged.
In contrast, some properties studied here are intrinsic to $\coneJ$.
For example, Slater’s condition $\coneJ \cap  \SymMat^n_{++} \not= \emptyset$ % (J ∩ Sⁿ₊₊ ≠ ∅) 
and the minimal face of $\SymMat^n_+$ %Sⁿ₊
containing $\coneJ$ depend only on $\coneJ$ and not on a particular representation.
Redundant constraints, however, are features of a chosen $\BC$ rather than
structural properties of $\coneJ$.
Thus, our results provide verifiable sufficient conditions,
expressed in terms of a representation $\BC$,
that guarantee the intrinsic geometric property $\coneJ \in \wFC(\coneG^n)$, %J ∈ Fb(Γₙ),
{\it i.e.}, that $\coneJ$ is ROG.

\medskip

We have presented two types of sufficient conditions 
under which a semi-infinite QCQP is equivalent to its SDP relaxation.
The first type, condition (B), extends  
the  assumptions in \cite[Theorem 4.1]{ARIMA2023} and \cite[Proposition 1]{ARGUE2023} 
for QCQPs with finitely many inequality constraints 
to the case of infinitely many inequality constraints. 
Even for QCQPs with finitely many constraints, 
condition (B) is weaker than the corresponding conditions in those works.
Condition (B) may be viewed as a sufficient structural condition on a representation of $\coneJ$ 
ensuring that $\coneJ$ (or equivalently, its minimal face) is ROG.
The second type, condition (C), is a new variant of 
NIQCC. We have shown that (C) implies (B) in general, 
and that the two are equivalent under appropriate additional assumptions. 
Furthermore, we have  generalized condition (III) (the NIQCC from 
 \cite[Corollary 2]{JOYCE2024}) to condition (D). 
 
 \medskip

In addition to the equivalence results for semi-infinite QCQPs and their SDP relaxations, 
which have been studied here based on ROG  % \cite{ARGUE2023,BLEKHERMAN2017,HILDEBRAND2016} 
and NIQCC, % \cite{BURER2015,JOYCE2024,YANG2018}, 
there exist several other well-studied classes of 
QCQPs whose % convex relaxations, including 
SDP % , DNN, and CPP 
relaxations are exact. Notable examples include convex QCQPs, where both the objective and all quadratic constraints are convex, 
and QCQPs defined by specific sign pattern conditions on the objective and constraint matrices
\cite{AZUMA2022,AZUMA2023,KIM2003,SOJOUDI2014}. 
%Thus 
These QCQPs have a fundamentally different nature, and most existing studies on them are largely independent of 
the ROG and NIQCC properties investigated in this paper. 
In our recent paper~\cite{KOJIMA2025}, we have discussed 
how homogenized and non-homogenized NIQCC properties can be incorporated into such QCQPs.

\bigskip

\smallskip
 
 \noindent
 {\bf \Large Acknowledgments.}
 The authors gratefully acknowledge the Associate Editor and two anonymous referees for their constructive feedback and insightful comments, 
 which have significantly improved this paper.
 In particular, their suggestions drew our attention to important references 
 %their suggestions directed us to important literature 
 on non-homogenized non-intersecting quadratic constraints.

\bibliography{./enhFOM}

\end{document}